# HIGHER FUNDAMENTAL FUNCTORS FOR SIMPLICIAL SETS

*by Marco GRANDIS* (*)

**Résumé.** On introduit une théorie d'homotopie combinatoire pour le topos des *ensembles simpliciaux symétriques* (préfaisceaux sur les cardinaux finis positifs), en étendant une théorie développée pour les complexes simpliciaux [11]; comme avantage essentiel de cette extension, le groupoïde fondamental devient l'*adjoint à gauche* d'un foncteur nerf symétrique et *préserve les colimites*, une propriété forte de van Kampen. On a des résultats analogues en toute dimension $\leq \omega$.

On développe aussi une notion d'*homotopie orientée* pour les ensembles simpliciaux ordinaires, avec un foncteur n-*catégorie fondamentale* adjoint à gauche du n-*nerf*. Des constructions similaires peuvent être données dans plusieurs catégories de préfaisceaux.

## Introduction

Higher homotopy theory and higher category theory present well known formal similarities; interestingly, coherence problems appeared at the same time (1963), independently, in both domains (Stasheff [25, 26], Mac Lane [20]). To develop their links is thus a natural problem with a variety of facets, also investigated in Grothendieck's manuscript [13]. The present work is a part of this design: we want to obtain higher homotopy groupoids as left adjoint functors, which yields a strong, simple version of the Seifert - van Kampen theorem: *such functors preserve all colimits*. This cannot be achieved within topological spaces, nor – on the other hand – within simplicial complexes (cf. 2.6); but can be obtained in suitable presheaf categories: a sort of equilibrium point, representing well-behaved spaces in a sufficiently flexible categorical structure.

The cases we will consider here are the well known topos **Smp** = $\mathbf{Set}^{\mathbb{A}^{\mathrm{op}}}$ of *simplicial sets*, or presheaves X: $\mathbb{A}^{\mathrm{op}} \to$ **Set** on the category $\mathbb{A}$ of finite positive

(*) Work supported by MURST Research Projects.

*ordinals* (with monotone maps) and the topos $!\mathbf{Smp} = \mathbf{Set}^{!\mathbb{A}^{op}}$ of *symmetric simplicial sets*, or presheaves on the category $!\mathbb{A}$ of finite positive *cardinals* (with all mappings); characterisations of this site, similar to the classical ones for $\mathbb{A}$, can be found in [12], Section 2. The *cubical* and *globular* analogues are only mentioned (5.4). All these are instances of 'combinatorial toposes', a notion which it would be interesting to develop in full generality; some desirable features are abstracted in 5.4.

First, we extend to the presheaf category **!Smp** the intrinsic homotopy theory developed in [11] for the well-known category **Cs** of *simplicial complexes*, or *combinatorial spaces* (sets with distinguished finite parts); the latter can be identified with those symmetric simplicial sets where each entry is determined by its vertices (the *simple presheaves*, 1.3-4), forming a cartesian closed, full subcategory of **!Smp**. Homotopies in **!Smp** are still parametrised on the *combinatorial integral line* **Z**, with combinatorial structure defined by the *contiguity* relation $i!j$, $|i-j| \leq 1$ (2.1). The fundamental groupoid $\Pi_1$: **!Smp** $\to$ **Gpd** (2.5) is now left adjoint to a natural functor $M_1$: **Gpd** $\to$ **!Smp**, the *symmetric nerve* of a groupoid, and preserves all colimits; this general fact also contains (or explains) the standard version of van Kampen in **Cs** [11], since those *particular* pushouts are preserved by the embedding in **!Smp** (2.6).

Higher fundamental groupoids are dealt with in Section 3 and shown to be left adjoints (Pasting Theorem, 3.9). Their construction, rather complicated, exploits some lifting techniques for endofunctors of presheaves, which reduce their 'calculation' to the level of the underlying set-valued functors (1.5-6) and might be useful elsewhere. In the literature, various notions of weak or strict n-groupoids have been considered in connections with homotopy (cf. [3, 5, 14, 15, 28] and their references). Here, we prefer to use the strict notion (but a semi-strict version is sketched in 3.8).

Similarly, within **Smp**, simple presheaves form a cartesian closed subcategory $\uparrow$**Cs**, of *directed simplicial complexes* (sets with distinguished words). In Sections 4, 5 we introduce a notion of 'directed homotopy' in these two categories; its non-reversible homotopies are now parametrised on the *directed integral line* $\uparrow$**Z**, with simplicial structure produced by *consecutivity*, a reflexive relation $i \prec j$ ($i \leq j \leq i+1$). Again, the *fundamental category* functor $\uparrow\Pi_1$: **Smp** $\to$ **Cat** is left adjoint to the usual nerve $N_1$: **Cat** $\to$ **Smp** and preserves colimits. Applications to image analysis, in the same line as those of [11], are briefly mentioned (4.5). The main formal results concern higher fundamental n-categories; in the topos **Smp**, all these functors are left adjoints (5.1-3): the functor of higher reduced paths $\uparrow\Pi_\omega$: **Smp** $\to$ $\omega$-**Cat** has a right adjoint $\omega$-*nerve* $N_\omega$: $\omega$-**Cat** $\to$ **Smp**, whereas $\uparrow\Pi_n$ is left adjoint to the restricted nerve $N_n$: n-**Cat** → **Smp**. As discussed in 5.3, it would be interesting to determine whether $N_\omega$ coincides with Street's ω-nerve [27].

General comments on defining 'combinatorial homotopy' in presheaf categories can be found at the beginning (1.2) and the conclusion (5.4), and could be read now, as a sequel to this Introduction. The basic facts used here appear to be: (A) a finite product of representable presheaves is a finite colimit of representables; (B) the site contains a 'standard interval' and all representable objects are connected with respect to it. The condition (A), which appears everywhere in Algebraic Topology (a product of tetrahedra is a pasting of tetrahedra on lower tetrahedra), intervenes here formally to make 'special endofunctors' of presheaves 'work' (1.6).

Finally, let us remark that various notions appear here in two versions, the *directed* case versus the *symmetric* (or reversible) one, as in the opposition *category-groupoid*; such opposition will often be marked by a prefix, ↑ (or 'd-') versus ! (or 'r-'), affecting the (historically) secondary notion. This choice will not be homogeneous, since category theory gives priority to the directed notions while homotopy has been developed in a reversible way. Thus, the well known category **Smp** of simplicial sets is a *directed* notion (based on ordinals), opposed to a less used category !**Smp** of *symmetric* simplicial sets (based on cardinals); on the other hand, the usual fundamental groupoid $\Pi_1 X$ is a *reversible* notion, here opposed to a directed notion of fundamental *category* $\uparrow\Pi_1 X$ (of some directed structure). The term 'directed' is preferred to 'ordered' or 'oriented', which might be misleading (cf. 4.2).

**Acknowledgements.** In November 1998, at a PSSL meeting in Trieste, Bill Lawvere suggested I might extend my study of the homotopy of simplicial complexes to symmetric simplicial sets, on the basis of his draft [19] where the fundamental groupoid of the latter is presented as a left adjoint. I would like to express my gratitude for his kind encouragement and for helpful discussions.

**Notation**. Finite ordinals and finite cardinals coincide. The term 'graph' stands for *reflexive oriented graph*. The symbol $|-|$ denotes an underlying structure or the module of a real number, never a geometric realisation; $[i, j]$ is always an interval of the *integral* line **Z**; the index $\kappa$ varies in the set $\{-, +\}$.

## 1. Presheaves and simple presheaves

Within the category **!Smp** of symmetric simplicial sets, the *simple* presheaves (1.3) form the cartesian closed subcategory of simplicial complexes (1.4). Endofunctors of presheaves are studied in 1.5-6; this part can be omitted during a first reading.

**1.1. Simplicial toposes.** There are two main 'combinatorial sites' *of simplicial type*, the *directed* and the *symmetric* one.

(a) The well known *simplicial site* $\mathbb{\Delta}$ consists of all positive finite ordinals $[n] = n{+}1$ $(n \geq 0)$ and *monotone* mappings; it gives the topos $\mathbf{Smp} = \mathbf{Set}^{\mathbb{\Delta}^{op}}$ of (directed) simplicial sets. A presheaf $X = ((X_n), (\lambda^*))\colon \mathbb{\Delta}^{op} \to \mathbf{Set}$ consists of a sequence of sets $X_n = X[n]$, with actions $\lambda^*\colon X_m \to X_n$ for all *monotone* mappings $\lambda\colon [n] \to [m]$ (generated by faces and degeneracies, under the well known simplicial relations). The representable presheaf $\mathbb{\Delta}(-, [n])$ is usually written $\Delta[n]$.

(b) The *symmetric simplicial site* will be the full subcategory $!\mathbb{\Delta} \subset \mathbf{Set}$ of finite non-zero cardinals (containing $\mathbb{\Delta}$, with the same objects); it produces the topos $!\mathbf{Smp} = \mathbf{Set}^{!\mathbb{\Delta}^{op}}$ of *symmetric simplicial sets* $X = ((X_n), (\lambda^*))\colon !\mathbb{\Delta}^{op} \to \mathbf{Set}$, with actions $\lambda^*\colon X_m \to X_n$ for *all* mappings $\lambda\colon [n] \to [m]$; again, such actions can be described by generators and relations ([12], 2.3). In particular, the symmetric group $S_{n+1}$ operates on the right on $X_n$, by $x\rho = \rho^*(x)$.

(c) By n-truncation, the full subcategories $\mathbb{\Delta}_n$ and $!\mathbb{\Delta}_n$, with objects [0], [1],... [n], give the presheaf category $\mathbf{Smp}_n$, of n-*truncated simplicial sets*, and the symmetric analogue $!\mathbf{Smp}_n$.

**1.2. Discrete sites and intervals.** Even if the simplicial sites will be our main concern, some points are better viewed in a more general context. First, let $\mathcal{C}$ be a *discrete site*, in the sense that:

$(A_0)$ $\mathcal{C}$ *is a small category with a terminal object* $\bullet$ (also written [0]),

so that the presheaf category $\mathrm{Psh}\mathcal{C} = \mathbf{Set}^{\mathcal{C}^{op}}$ is made concrete by the forgetful functor

(1) $\quad U\colon \mathrm{Psh}\mathcal{C} \to \mathbf{Set}$,

$\quad U(X) = X[0] = X_\bullet = |X| = \mathrm{Lim}(X\colon \mathcal{C}^{op} \to \mathbf{Set})$.

To fix notation, the Yoneda embedding is written $y\colon \mathcal{C} \to \mathrm{Psh}\mathcal{C}$ and a representable presheaf as $y(c) = \hat{c} = \mathcal{C}(-, c)\colon \mathcal{C}^{op} \to \mathbf{Set}$; by Yoneda, $X(c) = (\hat{c}, X) = U(X^{\hat{c}})$. A *point*, or *vertex*, of a presheaf $X\colon \mathcal{C}^{op} \to \mathbf{Set}$ is an element of $|X| = X[0]$; the *vertices*, or 0-*dimensional faces*, of $x \in X(c)$ are the points $u^*(x)$, for all arrows $u: [0] \to c$ in $\mathcal{C}$, if any. In particular, the points of $\hat{c}$ are precisely the maps $u: [0] \to c$, and coincide with the vertices $u^*(\varepsilon_c)$ of its free generator, $\varepsilon_c = id(c)$.

The forgetful functor U has a chain of adjunctions (as considered in [19] for general toposes)

(2) $E \dashv D \dashv U \dashv C$, $\qquad$ $(DS)(c) = S$,

$\quad$ $(CS)(c) = S^{\mathcal{C}([0], c)}$, $\qquad$ $E(X) = \mathrm{Colim}(X: \mathcal{C}^{op} \to \mathbf{Set})$,

where $D: \mathbf{Set} \to \mathrm{Psh}\mathcal{C}$ sends a set S to its *discrete* (or constant, or diagonal) presheaf, while CS is the *codiscrete* (or chaotic) presheaf on S. The left adjoint E acquires a homotopical meaning, *and will be written as* $\Pi_0$, under the following condition on the site:

($B_0$) $\mathcal{C}$ *contains a 'standard interval'*, written [1], *linked to the terminal* [0] *by two* faces $\partial^-, \partial^+ : [0] \rightrightarrows [1]$; *each object* c *in* $\mathcal{C}$ *is* connected *with respect to this interval.*

The last condition will mean that c has some point $u: [0] \to c$ and any two of them are related, under the equivalence relation generated by uRu': there is some map $a: [1] \to c$ with $a\partial^- = u$, $a\partial^+ = u'$. Then, the set $\Pi_0 X = \mathrm{Colim}(X)$ is just the coequaliser of the faces $X[1] \rightrightarrows X[0]$.

Say that a presheaf X is *connected* if it satisfies the same condition above; or equivalently, if $\Pi_0 X$ is a singleton. Combinatorial aspects, different from the behaviour of topological spaces, already appear at this 0-dimensional level: $\Pi_0$ preserves finite products and all colimits, but does not preserve infinite products, generally. A presheaf decomposes as the sum of its connected components

(3) $X = \Sigma_\alpha X_\alpha$, $\qquad$ $X_\alpha(c) = \eta_c^{-1}\{\alpha\}$ $\qquad$ $(\alpha \in \Pi_0 X)$.

Moreover, the elementary interval produces the presheaf endofunctor $X^{[1]}$ of *elementary paths*, or immediate paths. Since these cannot be concatenated (in our examples), we construct in $\mathrm{Psh}\mathcal{C}$ each *finite interval* [i, j] (for integers $i \le j$) pasting $j - i + 1$ copies of [1]'s on $j - i$ copies of [0]'s; precisely, as a *finite, connected colimit of representables*, as exemplified below for [2, 4]

(4)

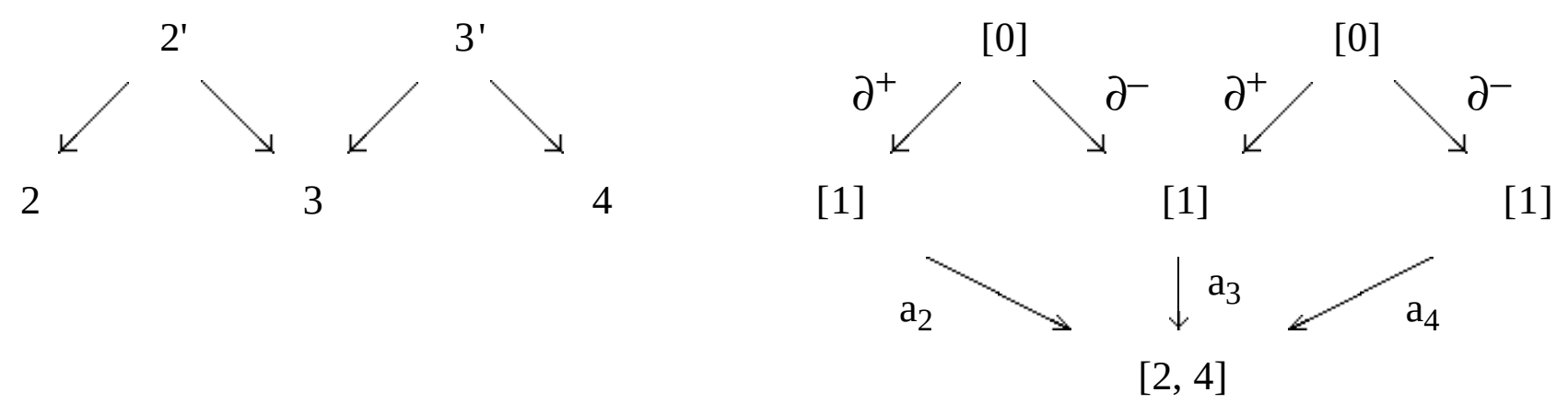

**N** and **Z** are constructed as similar countable colimits, and the representable presheaves y[0], y[1] are identified with [0, 0] and [0, 1], respectively. The path endofunctor can now be defined as a suitable colimit of endofunctors $X^{[i, j]}$ (cf. Sections 3, 5; such colimits are studied in 1.5-6).

We shall not proceed in this abstract way. In the simplicial toposes, it will be much simpler to use a *direct description* of the intervals rather than this colimit presentation; however, the fact that all finite intervals are finite colimits of representables will be crucial. The condition ($B_0$) holds in all simplicial (non-augmented!) cases, including the truncated ones. It is *not* assumed below (while ($A_0$) is), but will be reconsidered at the end (5.4), with a richer algebraic structure on the interval.

**1.3. Simple presheaves.** The presheaf X will be said to be *simple* if its elements are determined by their (indexed) family of vertices: if $x, x' \in X(c)$ and $u^*(x) = u^*(x')$ for all $u: [0] \to c$, then $x = x'$. The full subcategory of simple presheaves $SPsh\mathcal{C} \subset Psh\mathcal{C}$ has a reflector (left adjoint to the embedding)

(1) $p: \mathbf{Set}^{\mathcal{C}op} \to SPsh\mathcal{C}, \qquad (pX)(c) = (X(c))/{\sim_c},$

obtained by identifying, in each component, the elements which have the same family of vertices; plainly, p preserves finite products. Therefore, $SPsh\mathcal{C}$ is *cc-embedded* in the presheaf category $\mathbf{Set}^{\mathcal{C}op}$ and *cartesian closed 'with the same exponentials'*. (Cc-embedding means here a reflective full subcategory $j: \mathbf{A} \to \mathbf{X}$ of a cartesian closed category **X**, whose reflector $p: \mathbf{X} \to \mathbf{A}$ *preserves finite products*; then, by a Freyd's theorem ([7], 1.31), also **A** is cartesian closed, with the 'same' exponentials of **X**, $A^B = p((jA)^{jB})$, preserved by j.) In the simplicial cases (directed or symmetric), all representable presheaves are simple.

(A more general notion has been considered by Lawvere [19]. In a concrete topos $U: \mathcal{X} \to \mathbf{Set}$, with left adjoint $D \dashv U$, an object X is said *to live on its points* if each counit-map $DUY \to Y$ yields an *injective* mapping $\mathcal{X}(Y, X) \to \mathcal{X}(DUY, X)$. For $\mathcal{X} = \mathbf{Set}^{\mathcal{C}op}$, it is easy to see that this is equivalent to X being a simple presheaf.)

**1.4. Simplicial complexes.** We give now a direct description of the cc-embedded subcategories of simple presheaves in the toposes !**Smp**, !$\mathbf{Smp}_1$ (full, reflective, cartesian closed, with the same products and exponentials; 1.3). The directed analogue will be treated at the beginning of Section 4.

First, a *simplicial complex*, or *combinatorial space*, is a set X equipped with a set !X of finite subsets of X, called *linked parts* (or simplices), which contains the empty subset and all singletons, and is down-closed: if $\xi$ is linked, any $\xi' \subset \xi$ is so. A *morphism* of combinatorial spaces, or *map*, or *combinatorial mapping* f: $X \to Y$ is a mapping which preserves the linked parts.

The category **Cs** of combinatorial spaces is cartesian closed, with all limits and colimits (computed as in **Set** and equipped with the initial or final structure determined by the structural maps, respectively); it is regular, but not Barr-exact; the 2-point codiscrete object $[0, 1] = C\{0, 1\}$ classifies the *regular subobjects* (subsets equipped with the induced structure, the initial one). The internal hom, $\mathrm{Hom}(X, Y) = Y^X$, is the set of morphisms **Cs**(X, Y), where a finite subset of maps $\varphi$ is linked if, for all $\xi$ linked in X, $\varphi(\xi) = \cup_{f\in\varphi} f(\xi)$ is linked in Y.

The realisation C: $!\mathbb{\Delta} \to$ **Cs** of cardinals as *codiscrete* combinatorial spaces (all subsets are linked) yields a canonical cc-embedding

(1) c: **Cs** $\to$ !**Smp**,

$(cA)_n$ = !**Smp**(C[n], A) = {[n]-words in some linked part},

which identifies a combinatorial space with a *simple* symmetric simplicial set. Its reflector is

(2) u: !**Smp** $\to$ **Cs**, $\qquad$ $uX = X_0$,

where a non-empty part $\xi$ is linked in $X_0$ iff $\xi = \{i^*x \mid i: [0] \to [n]\}$, for some $x \in X_n$.

More simply, a *tolerance set* is a set equipped with a tolerance relation x!x' (reflexive and symmetric); they form the category **Tol** (with obvious mappings), which again is complete, cocomplete and cartesian closed; the exponential $Y^X$ has f!g whenever x!x' implies f(x)!g(x'). Plainly, a tolerance set is 'the same' as a simple involutive graph (simple in the usual sense of graph theory), i.e. a simple presheaf in !$\mathbf{Smp}_1$.

We have thus a commutative square of cc-embeddings

$$(3) \qquad \begin{array}{ccc} \mathbf{Tol} & \underset{c}{\overset{t}{\leftrightarrows}} & \mathbf{Cs} \\ c \downarrow\uparrow u & & c \downarrow\uparrow u \\ !\mathbf{Smp}_1 & \underset{cosk}{\overset{tr}{\leftrightarrows}} & !\mathbf{Smp} \end{array}$$

where the embedding c: **Tol** $\to$ **Cs** equips a tolerance set with the combinatorial structure making a finite part linked iff all its pairs belong to the tolerance relation; its left adjoint t: **Cs** $\to$ **Tol** is the forgetful functor: x!x' iff $\{x, x'\} \in$ !X.

The representable presheaf $!\Delta[n]$ is identified with the codiscrete object $C[n]$ on the set $[n] = \{0,... n\}$, in each of these categories: the codiscrete tolerance set (all pairs are linked), the codiscrete combinatorial space (all subsets are linked), the codiscrete involutive graph.

**1.5. Endofunctors of presheaves.** Coming back to the general situation, the discrete site $\mathcal{C}$ is again a small category with terminal object $[0] = \bullet$. The core of our construction will consist of endofunctors $F\colon \mathrm{Psh}\mathcal{C} \to \mathrm{Psh}\mathcal{C}$ (path functors, their powers and other derived operators), and we need to reduce them to their 'set-valued level' $UF\colon \mathrm{Psh}\mathcal{C} \to \mathbf{Set}$. The end of this Section is devoted to this point (and can be omitted, now; further use of this part will always be mentioned).

Let us begin noting that, in the (very large) categories of functors $(\mathrm{Psh}\mathcal{C}, \mathbf{Set})$ and $(\mathrm{Psh}\mathcal{C}, \mathrm{Psh}\mathcal{C})$, (small) limits and colimits exist and are pointwise calculated (in $\mathbf{Set}$ and $\mathrm{Psh}\mathcal{C}$; whence in $\mathbf{Set}$ again). Therefore, in both cases, finite limits commute with (small) filtered colimits, while *connected* finite limits commute with *pseudo*-filtered colimits, or – equivalently – with sums and filtered colimits [21].

The following well known procedure will be of much use: the *canonical lifting* of a set-valued functor $F_\bullet\colon \mathrm{Psh}\mathcal{C} \to \mathbf{Set}$ to an endofunctor $F = \mathcal{L}(F_\bullet)$ such that $UF = F_\bullet$

(1) $F\colon \mathrm{Psh}\mathcal{C} \to \mathrm{Psh}\mathcal{C}, \qquad\qquad (FX)(c) = F_\bullet(X^{\hat{c}}),$

corresponding, more formally, to the following composite $F'\colon \mathcal{C}^{\mathrm{op}} \times \mathrm{Psh}\mathcal{C} \to \mathbf{Set}$

(2) $$\mathcal{C}^{\mathrm{op}} \times \mathrm{Psh}\mathcal{C} \xrightarrow{\;y^{\mathrm{op}} \times 1\;} (\mathrm{Psh}\mathcal{C})^{\mathrm{op}} \times \mathrm{Psh}\mathcal{C} \xrightarrow{\;\mathrm{Hom}\;} \mathrm{Psh}\mathcal{C} \xrightarrow{\;F_\bullet\;} \mathbf{Set}\,.$$

Any natural transformation $\alpha_\bullet\colon F_\bullet \to G_\bullet\colon \mathrm{Psh}\mathcal{C} \to \mathbf{Set}$ has a canonical lifting $\alpha\colon F \to G$, with $(\alpha X)(c) = \alpha_\bullet(X^{\hat{c}})$, which is componentwise monic (resp. epi) if and only if $\alpha_\bullet$ is so. This gives a retraction of our functor categories

(3) $\mathcal{L} : (\mathrm{Psh}\mathcal{C}, \mathbf{Set}) \rightleftarrows (\mathrm{Psh}\mathcal{C}, \mathrm{Psh}\mathcal{C}) : \mathcal{U} = U\circ - \qquad\qquad (\mathcal{U}\mathcal{L} = 1);$

both $\mathcal{L}$ and $\mathcal{U}$ preserve (small) limits and colimits, since these are calculated pointwise. Thus, the endofunctors $F\colon \mathrm{Psh}\mathcal{C} \to \mathrm{Psh}\mathcal{C}$ which are canonical liftings *are closed under limits and colimits* in the functor category $(\mathrm{Psh}\mathcal{C}, \mathrm{Psh}\mathcal{C})$; they contain all representable endofunctors $(-)^Y\colon \mathrm{Psh}\mathcal{C} \to \mathrm{Psh}\mathcal{C}$. If $F_\bullet$ preserves all I-limits $X = \mathrm{Lim}\, X_i$ (based on a small category I), also $F$ does $(FX(c) = F_\bullet(X^{\hat{c}}) = \mathrm{Lim}\, F_\bullet(X_i{}^{\hat{c}}) = \mathrm{Lim}\,((FX_i)(c))$; the converse is trivial, and also holds for colimits.

A non-trivial problem is to determine a subclass of such endofunctors *closed under composition*, sufficiently large to include path functors of simplicial sets and their 'derived' constructs. Let us start noting that a functor $F: \mathrm{Psh}\mathcal{C} \to \mathrm{Psh}\mathcal{C}$ which preserves 'exponentials by representables'

(4) $F(X^{\hat{c}}) = (FX)^{\hat{c}}$, for all $c$ in $\mathcal{C}$,

is necessarily a canonical lifting. The condition (4) is stable under *composition* and *limits* in the functor category $(\mathrm{Psh}\mathcal{C}, \mathrm{Psh}\mathcal{C})$; it is satisfied by all representable endofunctors $(-)^Y$.

If the discrete site $\mathcal{C}$ has finite products (which is true for $!\mathbb{A}$ but not for $\mathbb{A}$), then this condition is also *necessary* for canonical liftings, which are therefore closed under composition

(5) $(F(X^{\hat{c}}))(d) = F_\bullet(X^{\hat{c}\times\hat{d}}) = (FX)(c\times d) = |(FX)^{\hat{c}\times\hat{d}}| = ((FX)^{\hat{c}})(d)$.

**1.6. Special endofunctors.** We consider now a weaker condition on the site and a restricted class of liftings. It is well known that a presheaf $A$ is *finitely presentable* (i.e. the hom-functor $(A, -): \mathrm{Psh}\mathcal{C} \to \mathbf{Set}$ preserves filtered colimits [8, 23, 2] if and only if it is a finite colimit of representables. (In the simplicial cases, this means finitely generated.) Similarly, we say that $A$ is *connected finitely presentable* if it is a connected finite colimit of representables, or equivalently if $(A, -)$ preserves pseudo-filtered colimits.

A condition of interest on our discrete site, obviously satisfied by the simplicial ones and many others (see comments in 5.4), will be that, equivalently:

(A) *each product* $\hat{c}\times\hat{d}$ *of representable presheaves in* $\mathrm{Psh}\mathcal{C}$ *is connected finitely presentable,*

(A') *connected finitely presentable objects in* $\mathrm{Psh}\mathcal{C}$ *are closed under finite products,*

(A") *if* $A$ *is connected finitely presentable, the functor* $(-)^A: \mathrm{Psh}\mathcal{C} \to \mathrm{Psh}\mathcal{C}$ *preserves pseudo-filtered colimits,* where the first equivalence is obvious and the second comes from $(X^A)(c) = \mathrm{Psh}\mathcal{C}(A\times\hat{c}, X)$. These conditions imply that finitely presentable objects are closed under finite products; the converse holds in the presence of $(B_0)$, as it follows easily from the preservation properties of $\Pi_0$ (1.2).

Let us assume that the discrete site $\mathcal{C}$ satisfies the condition (A) and say that $F: \mathrm{Psh}\mathcal{C} \to \mathrm{Psh}\mathcal{C}$ is a *special endofunctor* if it can be obtained as a *filtered colimit* of functors $F_iX = X^{A_i}$, with $A_i$ finitely presentable $(A: I \to (\mathrm{Psh}\mathcal{C})^{op}$, for $I$ a small filtered category); and *pseudo-special* if it can be obtained as a *pseudo*-filtered colimit of such functors, with $A_i$ *connected* finitely presentable.

The main properties we shall repeatedly use are as follows. The special endofunctors preserve exponentials by finitely presentable objects (stronger than 1.5.4) and are canonical liftings:

(1) $F(X^B) = \mathrm{Colim}_i\, (X^{A_i \times B}) = (\mathrm{Colim}_i\, X^{A_i})^B = (FX)^B$;

they always preserve finite limits; they are closed in $(\mathrm{Psh}\mathcal{C}, \mathrm{Psh}\mathcal{C})$ under filtered colimits and composition, which is computed by the following formula (for $F_iX = X^{A_i}$, $G_jX = X^{B_j}$)

(2) $GF(X) = \mathrm{Colim}_{I \times J}\, X^{A_i \times B_j}$ ;

the previous retraction (1.5.3) restricts to an isomorphism between special endofunctors and *special set-valued functors* $F_\bullet$, defined by the property $F_\bullet X = \mathrm{Colim}\, |X^{A_i}|$ (as above). Analogous facts hold for *pseudo*-special endofunctors, *connected* finitely presentable objects and *connected* finite limits.

Finally, let us note that – always assuming (A) – *within the endofunctors* $F\colon \mathrm{Psh}\mathcal{C} \to \mathrm{Psh}\mathcal{C}$ *which preserve connected finite limits*, the ones which are canonical liftings are still characterised by the fact of preserving exponentials by representables (1.5.4). The proof is as in 1.5.5, but the equality $F_\bullet(X^{\hat{c} \times \hat{d}}) = |F(X)^{\hat{c} \times \hat{d}}|$ derives now from (A), together with the fact that both functors $|F(X^{(-)})|$ and $|(FX)^{(-)}|$ take connected finite colimits to limits.

## 2. Combinatorial homotopy for symmetric simplicial sets, in degree 1

Combinatorial homotopy in **Cs** has been studied in [11]. Here we show how the path functor can be extended to symmetric simplicial sets. In this extension, the fundamental groupoid is left adjoint to a natural symmetric-nerve functor (2.6). Higher homotopy is deferred to the next section.

**2.1. The combinatorial integral line.** Consider again the commutative square of cc-embeddings linking the toposes **!Smp**, $\mathbf{!Smp}_1$ with their subcategories **Cs**, **Tol** of simple presheaves (1.4.3).

Following the outline in 1.2, combinatorial homotopy will be defined in each of these categories by means of the *combinatorial (integral) line* **Z**, the set of integers with the structure of *contiguity*, its finite intervals and their powers. Precisely, **Z** is a tolerance set, with tolerance relation $i!j$ if $|i - j| \leq 1$; as an involutive graph, it has precisely one map $i \to j$ whenever $i!j$; as a combinatorial space, a finite subset $\zeta \subset \mathbf{Z}$ is linked iff it is contained in some pair $\{i, i+1\}$; finally, as a symmetric simplicial set, its n-simplices are sequences $(i_0, \ldots i_n)$ of integers contained in some pair $\{i, i+1\}$. Thus, $\mathbf{Z} = sk_1tr_1\mathbf{Z}$ (where $tr_1$: $!\mathbf{Smp} \to !\mathbf{Smp}_1$ is the truncation, and $sk_1 \dashv tr_1$).

An integral interval $[i, j]$ has the induced structure. The *standard elementary interval* $[0, 1] = C\{0, 1\} \subset \mathbf{Z}$ is the codiscrete combinatorial space on two points, i.e. the representable presheaf $!\Delta[1]$. Of course, each *finite* interval $[i, j]$ is a finitely presentable presheaf, as a finite pasting of $\Delta[1]$'s on $\Delta[0]$'s.

For each $k \geq 1$, there is in $!\mathbf{Smp}$ and $!\mathbf{Smp}_1$ a k-*point combinatorial circle* $C_k$, the coequaliser of the faces of $[0, k]$; it only belongs to **Cs** for $k \geq 3$, to **Tol** for $k \geq 4$. On the other hand, *there is no standard circle*, i.e. no object representing the fundamental group functor in the appropriate homotopy category: either of pointed combinatorial spaces, as proved in [11], 6.5, or of pointed symmetric simplicial sets; in fact, such an object should have a 'generating loop' of a given finite length, which could not 'cover' longer loops, in some $C_k$. *This seems to be another distinctive feature of combinatorial homotopy*. Higher spheres will be considered in 2.7.

**2.2. Lines, paths, homotopies.** In the general case of symmetric simplicial sets, we have the representable endofunctor of *lines* (a canonical lifting of its set-valued part, 1.5)

(1) $L = (-)^{\mathbf{Z}}: !\mathbf{Smp} \to !\mathbf{Smp}, \qquad (LX)_n = (X^{\mathbf{Z}})_n = !\mathbf{Smp}(\mathbf{Z}, X^{!\Delta[n]})$.

In its underlying set

(2) $|LX| = !\mathbf{Smp}(\mathbf{Z}, X) = !\mathbf{Smp}(sk_1tr_1\mathbf{Z}, X) = !\mathbf{Smp}_1(tr_1\mathbf{Z}, tr_1X)$,

a line $a: \mathbf{Z} \to X$ amounts to its 1-truncation $(a_i: \mathbf{Z}_i \to X_i)_{i\leq 1}$, a morphism of involutive graphs, or equivalently a sequence of *consecutive* 1-simplices $a_1(i, i+1) \in X_1$ $(\partial^+a_1(i-1, i) = \partial^-a_1(i, i+1))$; of course, if $X$ is a *simple* presheaf, it suffices to assign the 0-component $a_0$, i.e. a map $a: \mathbf{Z} \to X$ of simplicial complexes: a sequence of points $a(i, i+1)$ where $\{a(i), a(i+1)\} \in !X$.

This line is said to be a *path* if it is eventually constant at the left and the right: there exists a *support* $\rho = [\rho^-, \rho^+] \subset \mathbf{Z}$, i.e. a finite integral interval such that the 1-simplex $a_1(i, i+1)$ is degenerate for $i < \rho^-$ and for $i \geq \rho^+$; then $\rho(a)$ will denote the *standard support* of the path $a$, i.e. the least admissible one (for a constant path, we choose $[0, 0]$). Paths form a subfunctor $P_\bullet \subset L_\bullet$: $!\mathbf{Smp} \to \mathbf{Set}$, which lifts canonically (1.5) to a *path endofunctor*, subfunctor of $L$

(3) $P: !\mathbf{Smp} \to !\mathbf{Smp}, \qquad (PX)_n = P_\bullet(X^{!\Delta[n]})$.

The trivial path at $x \in |X|$ is written $0_x = e(x)$, while $-a = r(a)$ denotes the reverse path: $(-a)_1(i, i+1) = a_1(-i, -i-1) = r(a_1(-i-1, -i))$.

The line-functors $LX = X^{\mathbf{Z}}$ are defined in the same way for the cc-subcategories **Tol**, **Cs**, !**Smp**$_1$ of diagram 1.4.3, and are restrictions of the previous L (since cc-embeddings preserve the exponential). A path is always defined as a line which is a path of the symmetric simplicial set. For combinatorial spaces, we find the same notion of line and path used in [11]: a line $a: \mathbf{Z} \to X$ is a combinatorial mapping (preserving the linked parts), and a finite set $\alpha$ of lines is linked iff each subset $\{a(i), a(i+1) \mid a \in \alpha\}$ is linked in X. PX is the subobject of those lines $a: \mathbf{Z} \to X$ which are eventually constant at the left and the right.

A *homotopy* $a: f \to g: X \to Y$ *of symmetric simplicial sets* is a map $a: X \to PY$ with $\partial^- a = f$, $\partial^+ a = g$; it is *bounded* if it factorises through some $Y^{[i, j]} \subset PY$ and *immediate* if it factorises through $Y^{[0, 1]}$; one can give an analytic description of the last notion, similar to the usual one for ordinary simplicial sets [24]. It is also interesting to note that the integral line is not contractible with respect to *bounded* homotopies ([11], 3.4). The general study of Homotopy given in [11] for combinatorial spaces, including homotopy pullbacks and the fibre sequence, can likely be extended to !**Smp**, but this cannot be done here, where our main concern is about higher fundamental groupoids.

**2.3. Delays.** To define a good concatenation of paths we have to collapse their constant parts. Consider the submonoid $D_1 \subset$ **Cs**(**Z**, **Z**) = !**Smp**(**Z**, **Z**) *of delays* ([11], 2.6), generated by the family of *elementary delays* $\delta_i$ (at $i \in \mathbf{Z}$)

(1) $\delta_i(t) = t$ if $t \leq i$, $\delta_i(t) = t - 1$ otherwise;

in other words, a delay $d: \mathbf{Z} \to \mathbf{Z}$ is an increasing surjective mapping (automatically combinatorial) which is the identity for t sufficiently small and strictly increasing for t sufficiently large. The crucial point is a sort of 'simplicial identity'

(2) $\delta_i.\delta_{j+1} = \delta_j.\delta_i$ $(i \leq j)$,

and the following consequence, called *the main property of delays* [11], or *cofiltering property*

(3) for any two delays $d_1, d_2$, there are delays $e_1, e_2$ such that $d_1e_1 = d_2e_2$;

(in fact, this holds *within* the generators by the simplicial identity; the general case follows by a repeated application of the particular one.)

A delay is not a path in $\mathbf{Z}$; but if $a: \mathbf{Z} \to X$ is a path, so is any delayed line $ad: \mathbf{Z} \to X$. Two paths $a, b \in |PX| \subset !\mathbf{Smp}(\mathbf{Z}, X)$ are said to be *congruent*, $a \equiv b$, if there exist two delays $d, d': \mathbf{Z} \to \mathbf{Z}$ such that $ad = bd'$. Congruence is an equivalence relation, by the main property above. Note that *any path is congruent to its translations*, since $a\delta_i$ is the path $a$ 'delayed of one unit' (for $i \leq \rho^-(a)$). Therefore, we can equivalently use the larger submonoid $\overline{D}_1 \subset \mathbf{Cs}(\mathbf{Z}, \mathbf{Z})$ of *generalised delays*, consisting of all increasing combinatorial mapping which are strictly increasing for $t$ sufficiently small or sufficiently large (generated by $D_1$ and elementary translations). In any congruence class of paths, an 'essential representative' is a path which has no repetitions on its standard support; such a representative is unique up to translations of $\mathbf{Z}$.

**2.4. Concatenation.** The set $\mathbf{s}$ of *supports* consists of the finite integral intervals $[i, j]$, with $i \leq j$; it will be equipped with a structure of commutative involutive monoid, in additive notation

(1) $[i, j] + [i', j'] = [i+i', j+j'], \qquad 0_{\mathbf{s}} = [0, 0], \qquad -[i, j] = [-j, -i];$

(note that the opposite pair $(-i, -j)$ does not denote an interval).

We want now to concatenate two *consecutive* paths $a, b \in |PX|$ $(\partial^+a = \partial^-b)$, extending the procedure used in [11] for simplicial complexes. Also here, this cannot be done in a natural way; any choice of a pair $\rho, \sigma$ of admissible supports for $a, b$ yields an *admissible concatenation* $c$ of our paths, with 'pasting point' at $\rho^+ + \sigma^-$ and admissible support $\rho + \sigma = [\rho^- + \sigma^-, \rho^+ + \sigma^+]$

(2) $c_1(i, i+1) = a_1(i - \sigma^-, i+1 - \sigma^-), \qquad \text{for } i < \rho^+ + \sigma^-,$
$c_1(i, i+1) = b_1(i - \rho^+, i+1 - \rho^+), \qquad \text{for } i \geq \rho^+ + \sigma^-.$

However, *all admissible concatenations of two given paths* $a, b$ *are congruent*: let $c$ be one of them, derived from supports $\rho, \sigma$; then, varying $\rho^-$ or $\sigma^+$ has no effect on $c$, while increasing $\rho^+$ (resp. $\sigma^-$) of one unit yields a concatenation $c' = c\delta_i$ delayed of one unit at a suitable instant.

To get a well defined operation, the (standard) *concatenation* $a+b$ will be the one produced by the standard supports of $a$ and $b$. This is coherent with the structure of $\mathbf{s}$

(3) $\rho(0_X) = [0, 0], \qquad \rho(-a) = [-\rho^+(a), -\rho^-(a)], \qquad \rho(a+b) = \rho(a) + \rho(b),$

and, for each object $X$, the graph $|X| \leftrightarrows |PX|$ is an involutive category, in additive notation

(4) $0 + a = a = a + 0$, $(a + b) + c = a + (b + c)$,
$-0 = 0$, $-(-a) = a$, $-(a + b) = -b - a$.

The standard concatenation is *natural up to delays*: a map $f: X \to Y$ takes a path $a \in |PX|$ to a path $fa \in |PY|$, for which $\rho(a)$ is an *admissible* support. Therefore $f.(a + b)$, being computed on the basis of $\rho(a)$, $\rho(b)$, is just an *admissible* concatenation of $fa$, $fb$

(5) $f.(a + b) \equiv fa + fb$,

and we have a non-natural family of concatenation mappings, which *cannot be lifted to presheaves*

(6) $kX: |PX| \times_{|X|} |PX| \to |PX|$, $k(a, b) = a+b$.

**2.5. The fundamental groupoid.** On the other hand, this concatenation is sufficient to define a (functorial) 1-dimensional fundamental groupoid, as we show now. Let us start from the 2-truncated cubical set produced by $P$ (only faces are shown, and $\kappa = \pm$)

(1) $|X| \leftleftarrows |PX| \leftleftarrows\leftleftarrows |P^2X|$ $\partial_1^\kappa = P\partial^\kappa$, $\partial_2^\kappa = \partial^\kappa P$,

where $|P^2X|$ is the set of *double paths* $A: \mathbf{Z}^2 \to X$ (having a finite support $\rho\times\sigma$), with four faces $\partial_1^\kappa$, $\partial_2^\kappa$. This cubical set restricts, in the usual way, to a (truncated) globular set

(2) $|X| \leftleftarrows\!\rightarrow |PX| \leftleftarrows\!\rightarrow P_2(X)$

where $P_2(X) \subset |P^2(X)|$ is the subset of 2-*paths*, i.e. those double paths whose faces $\partial_1^\kappa = P\partial^\kappa$ are degenerate (the faces $\partial_2^\kappa = \partial^\kappa P$ being kept in the structure). As usual, two paths $a, b \in |PX|$ are said to be 2-*homotopic*, or *homotopic with fixed end-points*, if there exists a 2-homotopy $A: a \simeq_2 b$, i.e. some 2-path $A \in P_2(X)$ with $\partial^- P(A) = a$, $\partial^+ P(A) = b$; 2-homotopy is an equivalence relation, because $|P(PX)|$ is a category, with faces $\partial^\kappa P$.

Congruence up to delays implies 2-homotopy, whence the concatenation in $\Pi_1 X$ is natural. Indeed, given a path $a: \mathbf{Z} \to X$, there is a *caterpillar* 2-homotopy (as in [11], 2.8) $A: a\delta_i \simeq_2 a$, which modifies $a\delta_i$ taking the delay to the right of the support of $a$, where it is ineffective: let $A = aB_{ij}$, with $j \geq i \vee \rho^+(a)$ and $B_{ij}$: $\mathbf{Z}^2 \to \mathbf{Z}$, $B_{ij}(s, t) = \delta_{(i\vee t)\wedge j}(s)$

| | | | | | | | | | | |
|---|---|---|---|---|---|---|---|---|---|---|
| | … | i | i+1 | i+2 | … | j–1 | j | j+1 | … | ($s \in \mathbf{Z}$) |
| $\delta_j$: | … | i | i+1 | i+2 | … | j–1 | **j** | **j** | … | (t = j) |
| | … | i | i+1 | i+2 | … | **j–1** | **j–1** | j | … | |
| (3) | … | | | | | | | | … | |
| | … | i | **i+1** | **i+1** | … | j–2 | j–1 | j | … | (t = i+1) |
| $\delta_i$: | … | **i** | **i** | i+1 | … | j–2 | j–1 | j | … | (t = i). |

($B_{ij}$ is not a double path, but $aB_{ij}$ is a 2-path; the 'caterpillar wave' is shown in boldface letters.)

Finally, the coequaliser of the faces $P_2(X) \rightrightarrows |PX|$ yields the *fundamental groupoid* $\Pi_1 X = |PX|/\simeq_2$ of the symmetric simplicial set X. In fact, we know that it is an involutive category, functorial in X; we could prove directly, as in [11], 2.9, that reversed paths are inverses up to 2-homotopy; but this will follow from a general treatment of the fundamental n-groupoids $\Pi_n X$ (3.8).

**2.6. The adjunction.** It is easy to see that, for simplicial complexes, $\Pi_1$: **Cs** → **Gpd** does not preserve all pushouts and is *not* a left adjoint: pasting in **Cs** two copies of the elementary interval [0, 1] = !Δ[1] on their vertices just gives [0, 1], which is not consistent with the pasting of the corresponding fundamental groupoids. On the other hand, in the more flexible frame !**Smp** (or !$\mathbf{Smp}_1$), the same pasting gives the 'homotopically correct' answer, the 2-point circle $C_2$ (2.1). Again, the coequaliser of the faces $[0] \rightrightarrows [0, 1]$ is the point in **Cs**, but $C_1$ in !**Smp**.

And indeed, as noted in [19], symmetric simplicial sets do have an adjunction

(1) $\Pi_1$ : !**Smp** $\rightleftarrows$ **Gpd** : $M_1$ $\qquad\qquad$ $\Pi_1 \dashv M_1$.

Here, $M_1G$ is the *symmetric nerve* of a groupoid, i.e. the representable functor produced by the embedding c: !𝔸 → **Gpd** which sends the cardinal [n] to its codiscrete groupoid c[n]

(2) $M_1G = \mathbf{Gpd}(c(–), G)\colon\ !\mathbb{A}^{op} \to \mathbf{Set}$,

$(M_1G)_n = \mathbf{Gpd}(c[n], G) =$

$= \{(a_{ij})_{0\leq i,j\leq n} \mid a_{ij} \in G_1;\ a_{ii} \text{ is an identity};\ a_{ij} + a_{jk} = a_{ik}\}$,

$\lambda^*\colon (M_1G)_n \to (M_1G)_m, \qquad (\lambda^* a) = (a_{\lambda i, \lambda j})_{0\leq i,j\leq m} \qquad (\lambda\colon [m] \to [n])$.

As a strong and particularly simple version of the van Kampen theorem, it follows that here *the fundamental groupoid preserves all colimits*. This should be compared with a (standard) version of van Kampen in **Cs**, proved in [11], 6.4, and concerned with the preservation of *particular* pushouts (the object $X = U \cup V$ is the union of its subobjects U, V, hence their pushout over $U \cap V$); noting that such pushouts are preserved by the embedding in !**Smp**, we see that this version is a trivial consequence of the general result in !**Smp**, and 'explained' by it.

Similarly, but this already holds in **Cs** [11], the functor $\Pi_0$: !**Smp** → **Set** of path components is left adjoint to the discrete embedding D: **Set** → !**Smp** (as already remarked in 1.2).

**2.7. Combinatorial spheres.** Let $n \geq 1$. In **Cs**, various combinatorial models of the n-sphere can be given: see [11], 1.3, for a description of the *simplicial* n-*sphere* $\Delta S^n$, the *cubical* n-sphere $\square S^n$, and the *octahedral* n-sphere $\lozenge S^n$. More generally, one could consider a system $\square_k S^n \prec [0, k]^{n+1}$ of (pointed) cubical n-spheres providing $\pi_n X = \mathrm{Colim}_k[\square_k S^n, X]$; but the n-dimensional loop $p^n_k$: $\mathbf{Z}^n \to \square_k S^n$ required for this representation is somewhat complicated.

This purpose can be more easily achieved using a different model, the k-*collapsed* n-*sphere* $C_k S^n = \mathbf{Z}^n/\sim_k$, obtained by collapsing to a (base) point all the points of $\mathbf{Z}^n$ out of the cube $[1, k]^n$ (for $k \geq 2$). This object can be viewed as the surface of a 'pyramid' with basis $[1, k]^n$ and vertex the base-point [0]; its linked parts are the ones of the basis $[1, k]^n \subset \mathbf{Z}^n$, together with all 'triangles' $\xi \cup [0]$ where $\xi$ is linked in some face of the basis. The quotient map $p^n_k$: $\mathbf{Z}^n \to C_k S^n$ is the generator of $\pi_n(C_k S^n)$ (as it follows from the geometric realisation and [11], thm. 6.3); the system of quotient maps $q^n_k$: $C_{k+1} S^n \to C_k S^n$ $(k \geq 2)$ induced by the identity of $\mathbf{Z}^n$ *surrogates* a standard n-sphere, by an isomorphism $\pi_n X \cong \mathrm{Colim}_k[C_k S^n, X]$ induced by the quotient-maps $p^n_k$: $\mathbf{Z}^n \twoheadrightarrow C_k S^n$.

**2.8. Directed objects.** This homotopy theory can be lifted to ordinary simplicial sets X, by means of the (symmetric) intervals U[i, j] produced by the restriction functor U: !**Smp** → **Smp**; but U ⊣ sym (where sym(X) = Ran(X) is the right Kan extension of X: $\mathbb{\Delta}^{op}$ → **Set** along $\mathbb{\Delta}^{op} \subset\ !\mathbb{\Delta}^{op}$), and this lifting simply amounts to *combinatorial homotopy for the symmetric objects* sym(X): **Smp**(U[i, j], X) = !**Smp**([i, j], symX). Finer results will be obtained below (Sections 4, 5), by the *directed* standard interval ↑[0, 1] = Δ[1].

## 3. Higher fundamental groupoids for symmetric simplicial sets

Higher fundamental groupoids are introduced (3.7-8) and proved to be left adjoints (3.9). The lifting techniques for endofunctors of presheaves (1.5-6) are heavily used.

**3.1. Truncating maps.** To proceed to higher fundamental groupoids we need a more formal study of the path functor, with three new endofunctors of **!Smp** related to P by natural transformations

(1) $\mathbb{P} \to \mathrm{P} \to \mathbf{P} \to \mathcal{P}$;

all four can be obtained as a colimit of exponentials $X^{[i, j]}$, over four categories of supports $\mathbf{s} \subset \mathbf{t} \subset \mathbf{d} \subset \mathbf{c}$. The functor of Moore paths $\mathbb{P}$ will play an auxiliary role for the functors of reduced paths **P** and of strongly reduced paths $\mathcal{P}$; the latter will produce our higher fundamental groupoids, while **P** gives a sort of *semi-strict* n-groupoid (3.8). (The 'homotopies' that the new functors represent have various disadvantages with respect to P-homotopies, discussed for **Cs** in [11], 4.6.)

Let us begin reconsidering the ordinary path functor P. Recall that, for a symmetric simplicial set X, a path $a: \mathbf{Z} \to X$ amounts to a sequence of consecutive 1-simplices $a_1(i, i+1) \in X_1$ $(\partial^+a_1(i-1, i) = \partial^-a_1(i, i+1))$, which are degenerate when $i < \rho^-$ or $i \geq \rho^+$. This means that the set |PX| can be formalised as a colimit of sets $|X^{[i, j]}|$, with respect to trivial prolongations $|X^{[i, j]}| \to |X^{[i', j']}|$ of paths at end-points. To make this precise, consider the family of *truncating maps*, in **Tol** $\subset$ **!Smp**

(2) $p: \mathbf{Z} \to [i, j], \qquad p(t) = (i \vee t) \wedge j = i \vee (t \wedge j) \qquad (i \leq j),$

and the category $\mathbf{t}^{\mathrm{op}} \subset \mathbf{Tol} \subset \mathbf{!Smp}$ of their restrictions (written $\mathbf{t}^{\mathrm{op}}$ for convenience)

(3) $p: [i', j'] \to [i, j], \qquad p(t) = (i \vee t) \wedge j = i \vee (t \wedge j) \qquad (i' \leq i \leq j \leq j').$

The opposite **t** is a directed ordered set (whence a filtered category), corresponding to the inclusion relation $[i, j] \subset [i', j']$, with joins $([i', j'] \vee [i'', j''] = [i' \wedge i'', j' \vee j''])$; the maps (2) form a cone from **Z** to the inclusion functor $\mathbf{t}^{\mathrm{op}} \subset$ **!Smp**, and a cocone of $X^{(-)}: \mathbf{t} \to (\mathbf{!Smp}, \mathbf{!Smp})$.

Since a finite interval [i, j] is a finite colimit of representables (1.2.4), the path functor P and all its powers $\mathrm{P}^n$ are *special* endofunctors, determined by their set-valued part $|\mathrm{P}^n|$ (together with their natural transformations) and calculated by the composition formula 1.6.2

(4) $PX = \mathrm{Colim}_{\mathbf{t}} X^{[i, j]}$, $|PX| = \mathrm{Colim}_{\mathbf{t}}\, !\mathbf{Smp}([i, j], X)$,

(5) $P^n X = \mathrm{Colim}_{\mathbf{t}^n} X^{[i_1, j_1]\times ... \times[i_n, j_n]}$,

$|P^n X| = \mathrm{Colim}_{\mathbf{t}^n}\, !\mathbf{Smp}([i_1, j_1]\times ... \times[i_n, j_n], X)$.

The line functor $LX = X^{\mathbf{Z}}$ (which is not special, but is a canonical lifting) is linked to P by the natural transformation $PX \subset LX$ coming from the cone (2), and described at the 0-dimensional level by the inclusion $|PX| \subset |LX|$ of paths as 'lines with a finite support', considered from the beginning. Similarly, the derived natural transformation $P^n X \subset L^n X$ is described by a set-inclusion $|P^n X| \subset |L^n X|$ identifying an n-dimensional path $a \in |P^n X|$ with an n-dimensional line $a: \mathbf{Z}^n \to X$ having a finite support $[i_1, j_1]\times ... \times[i_n, j_n]$.

**3.2. The cubical comonad of paths.** The basic properties of homotopy repose on natural transformations linking the path endofunctor $P: !\mathbf{Smp} \to !\mathbf{Smp}$ and its powers. (A detailed description for simplicial complexes can be found in [11], Section 2).

First, the *degeneracy* $e: P \to 1$ and the two *faces* $\partial^\kappa: P \to 1$ $(\kappa = \pm)$ come from obvious natural transformations of the functors $X^{[i, j]}$

(1) $\partial^\kappa : [0] \rightrightarrows [i, j] : e$, $\partial^-(0) = i$, $\partial^+(0) = j$, $e(t) = 0$.

Deriving the remaining transformations from the colimit construction of P (3.1.4) would be rather long. We prefer to start from the structure of $\mathbf{Z}$ as an involutive lattice (without extrema!) *in* **Tol**

(2) $[0] \xleftarrow{e} \mathbf{Z} \overset{g^\kappa}{\leftleftarrows} \mathbf{Z}^2$ $r: \mathbf{Z} \to \mathbf{Z}$, $s: \mathbf{Z}^2 \to \mathbf{Z}^2$,

$g^-(t, t') = t \vee t'$, $g^+(t, t') = t \wedge t'$, $r(t) = -t$, $s(t, t') = (t', t)$,

deduce the corresponding transformations for the line functor $LX = X^{\mathbf{Z}}$ and restrict them to $P^n X \subset L^n X$; by canonical lifting (1.5), the fact that the restriction be legitimate need only be checked at the 0-dimensional level $|PX| \subset |LX|$ of paths, where it is obvious. We get thus the previous degeneracy e, two *connections* $g^-$, $g^+$, and two symmetries, the *reversion* r and the *interchange* s (of course, the faces $\partial^\kappa: P \to 1$ cannot be derived from $\mathbf{Z}$; one should replace it with the extended line $\overline{\mathbf{Z}}$ with two infinity points, a 'true' lattice). Globally, we have seven natural transformations

$$(3) \quad 1 \underset{e}{\overset{\partial^\kappa}{\leftleftarrows\!\rightarrow}} P \overset{g^\kappa}{\rightrightarrows} P^2 \qquad r: P \to P, \qquad s: P^2 \to P^2;$$

which satisfy the axioms of a *cubical comonad with symmetries* [9, 10] (the functorial dual-analogue of a *commutative, involutive cubical monoid*: a set equipped with two structures of commutative monoid $g^\kappa$, where the unit $\partial^\kappa$ of each operation is an absorbent element for the other, and the involution r turns each structure into the other):

(4) $\partial^\kappa.e = 1,\quad g^\kappa.e = Pe.e\ (= eP.e)$ (*degeneracy axiom*),
$Pg^\kappa.g^\kappa = g^\kappa P.g^\kappa \quad P\partial^\kappa.g^\kappa = 1 = \partial^\kappa P.g^\kappa$ (*associativity, unit*),
$P\partial^\kappa.g^{\kappa'} = e.\partial^\kappa = \partial^\kappa P.g^{\kappa'}$ (*absorbency*; $\kappa \neq \kappa'$),
$r.r = 1, \quad r.e = e, \quad \partial^-.r = \partial^+, \quad g^-.r = Pr.rP.g^+,$
$s.s = 1, \quad s.Pe = eP, \quad P\partial^\kappa.s = \partial^\kappa P, \quad s.g^\kappa = g^\kappa,$
$Pr.s = s.rP$ (*symmetries*).

These properties either come directly from the dual ones for **Z**, or (when faces are involved) are immediate at the level of set-valued functors. P preserves all *finite* limits (1.6), but does not preserve infinite products. In fact, this preservation only depends on |P| (1.5), and $(–)^{\mathbf{Z}}$ preserves all limits; now, a *finite* jointly monic family of maps $X \to X_i$ (e.g., the projections of a limit) reflects the lines which admit a support; but this obviously fails in the infinite case. Therefore, P has no left adjoint, and the homotopies which it generates have no cylinder functor.

**3.3. Moore paths.** Now, we want to modify P to obtain a natural concatenation. A first way is to consider paths with an assigned support, like Moore paths for topological spaces; the advantage on concatenation is paid with the existence of infinitely many constant paths at the same point.

Here, the endofunctor of *Moore paths*, or *paths with duration*, will be the following *pseudo-special* endofunctor: a sum over the set **s** of supports (2.4), i.e. a colimit on the discrete category **s**

(1) $\mathbb{P}X = \Sigma_{\mathbf{s}} X^{[i, j]}, \qquad |\mathbb{P}X| = \Sigma_{\mathbf{s}} !\mathbf{Smp}([i, j], X) \subset \mathbf{s}\times|PX|.$

Its set-valued part consists of all pairs $\hat{a} = (\rho, a)$ where $\rho \in \mathbf{s}$ is an admissible support of $a \in |PX|$. At each point $x \in |X|$ there are as many *constant paths* $(\rho, e(x))$ as supports $\rho \in \mathbf{s}$. A map $a: X \to \mathbb{P}Y$ amounts to a *homotopy bounded on each connected component* of X, with a choice of a support for each component. $\mathbb{P}$ has obvious faces, degeneracy, connections and symmetries, described below at set-level (adding, to the corresponding transformation of P, the obvious transformation of supports)

(2) $\partial^\kappa(\rho, a) = \partial^\kappa a$, $\quad e(x) = ([0, 0], ex)$, $\quad g^\kappa(\rho, a) = (\rho{\times}\rho, g^\kappa a)$,
$\quad r(\rho, a) = (-\rho, ra)$, $\quad s(\rho{\times}\sigma, a) = (\sigma{\times}\rho, sa)$.

Such transformations for $\mathbb{P}$ *nearly* form a cubical comonad with symmetries: the only axiom which fails (with respect to the list 3.2.4) is absorbency, 'because' of constant paths. On the other hand, we have a concatenation, defined on the pullback $\mathbb{Q}X = \mathbb{P}X \times_X \mathbb{P}X$ of consecutive pairs of Moore paths (still a pseudo-special functor) and determined by its 0-dimensional level

$$
(3) \qquad \begin{array}{ccc} \mathbb{Q}X & \xrightarrow{k^-} & \mathbb{P}X \\ {\scriptstyle k^+}\downarrow & & \downarrow{\scriptstyle \partial^+} \\ \mathbb{P}X & \xrightarrow{\partial^-} & X \end{array} \qquad \begin{array}{l} k\colon \mathbb{Q} \to \mathbb{P}\colon \mathbf{!Smp} \to \mathbf{!Smp}, \\ \\ k(\rho, a; \sigma, b) = (\rho+\sigma, a+b), \end{array}
$$

where $a+b$ is *now* the concatenation of $a$ and $b$ over the assigned supports $\rho, \sigma$ (2.4.2). The transformation $k$ does satisfy the usual axioms of concatenation [10]:

(4) $\partial^- k = \partial^- k^-$, $\quad \partial^+ k = \partial^+ k^+$, $\quad ke_{\mathbb{Q}} = e$,
$\quad kr_{\mathbb{Q}} = rk$, $\quad k\mathbb{P}.s' = s.\mathbb{P}k$,

where $e_{\mathbb{Q}}\colon X \to \mathbb{Q}X$, $r_{\mathbb{Q}}\colon \mathbb{Q}X \to \mathbb{Q}X$, and $s'\colon \mathbb{P}\mathbb{Q}X \to \mathbb{Q}\mathbb{P}X$ are induced by $e$, $r$, and $s$ respectively. Moreover, it is *semiregular*, i.e. it makes $|\mathbb{P}X|$ into a strict involutive category.

**3.4. Reduced paths.** We have already seen that a second way of getting a natural concatenation consists in identifying the paths of PX (or $\mathbb{P}$X) by the congruence up to delays, $\equiv$ (2.3-4).

A congruence class of paths $a^\bullet$ will be called a *reduced path* of X; we have thus a set-valued functor $\mathbf{P}_\bullet X = |PX|/\equiv$. Its lifting $\mathbf{P}\colon \mathbf{!Smp} \to \mathbf{!Smp}$, called the *reduced path functor*, can be obtained as a colimit based on the category $\mathbf{d} \subset \mathbf{Tol}^{op} \subset \mathbf{!Smp}^{op}$

(1) $\mathbf{P}X = \mathrm{Colim}_{\mathbf{d}}\, X^{[i,\, j]}$, $\quad |\mathbf{P}X| = \mathrm{Colim}_{\mathbf{d}}\, \mathbf{!Smp}([i, j], X) = |PX|/\equiv$,

where $\mathbf{d}^{op}$ consists of generalised delays $d\colon [i', j'] \to [i, j]$: all *increasing surjective mappings*, generated by translations and the (restricted) elementary delays $\delta_i\colon [h, k] \to [\delta_i(h), \delta_i(k)]$. *Unfortunately*, $\mathbf{d}$ is not filtered (all its arrows are monic, and $\mathbf{d}$ is not a preorder), but just 'diamond-filtered' and 'object-filtered' (any two arrows with the same domain and any two objects have a cocone): this follows easily from the fact that any $d: [i', j'] \to [i, j]$ can be decomposed as an 'ordinal' sum $d_i + ... + d_j$ of constant surjective mappings. And indeed **P** is *not* a pseudo-special endofunctor in the sense of 1.6 (as we prove in 3.6); it is nevertheless a canonical lifting.

The inclusion $\mathbf{t} \subset \mathbf{d}$ gives a natural projection (surjective at 0-level, whence at all levels)

(2) $h: PX \to \mathbf{P}X, \qquad h_\bullet: |PX| \to |\mathbf{P}X| = |PX|/\equiv.$

Since all admissible concatenations of two given paths a, b are congruent (2.4), the set-valued concatenation $|k|: |\mathbb{Q}| \to |\mathbb{P}|$ induces a natural transformation

(3) $k_\bullet: \mathbf{Q}_\bullet \to \mathbf{P}_\bullet, \qquad \mathbf{Q}_\bullet X = (|PX| \times_{|X|} |PX|)/(\equiv\times\equiv) = |\mathbf{P}X| \times_{|X|} |\mathbf{P}X|,$

which lifts to a concatenation, defined over the functor $\mathbf{Q}X = \mathbf{P}X \times_X \mathbf{P}X$ *of pairs of consecutive reduced paths* (1.5), so that $|\mathbf{P}X|$ is again a strict involutive category on $|X|$

(4) $$\begin{array}{ccc} \mathbf{Q}X & \xrightarrow{k^-} & \mathbf{P}X \\ {\scriptstyle k^+}\downarrow & & \downarrow{\scriptstyle \partial^+} \\ \mathbf{P}X & \xrightarrow{\partial^-} & X \end{array} \qquad\qquad \begin{array}{l} k: \mathbf{Q} \to \mathbf{P}: !\mathbf{Smp} \to !\mathbf{Smp}, \\ \\ \mathbf{Q}X = \mathbf{P}X \times_X \mathbf{P}X. \end{array}$$

**3.5. Strongly reduced paths.** One can transform the reversed 1-cells of $|\mathbf{P}X|$ in strict inverses, by a further identification of paths 'up to regressions', so that any concatenation $a - a$ be identified with the trivial path at $\partial^-a$ (as in the 2-groupoid for Hausdorff spaces of [14]).

To this effect, consider the semigroup $\Gamma_1 \subset \mathbf{Cs}(\mathbf{Z}, \mathbf{Z})$ of all *combinatorial surjective mappings* which are strictly increasing for t sufficiently small or large (possibly not increasing); it contains $\overline{D}_1$ (2.3) and is generated by translations, elementary delays and *elementary regressions* $\gamma_i: \mathbf{Z} \to \mathbf{Z}$

(1) $\gamma_i(t) = t \quad \text{if } t \leq i, \qquad \gamma_i(t) = t - 2 \quad \text{otherwise}.$

The latter satisfy the following identities

(2) $\gamma_i.\gamma_{j+2} = \gamma_j.\gamma_i \quad (i \leq j), \qquad \gamma_i.\gamma_{i+2} = \gamma_i.\gamma_{i+1} = \gamma_i.\gamma_i,$
$\delta_i.\gamma_{j+1} = \gamma_j.\delta_i \quad (i \leq j), \qquad \delta_i.\gamma_j = \gamma_j.\delta_{i+2} \quad (i \geq j),$

whence the semigroup $\Gamma_1$ satisfies the same cofiltering property of delays (2.3.2), and its action on paths produces a (stronger) congruence, defined similarly. Again, $\Gamma_1$-congruence implies 2-homotopy, because a regressed path $a\gamma_i = (... a_i, a_{i-1}, a_i,...)$ is immediately 2-homotopic to the delayed path $a\delta_i\delta_i = (... a_i, a_i, a_i,...)$.

More generally, the submonoid $\Gamma_n \subset \mathbf{Cs}(\mathbf{Z}^n, \mathbf{Z}^n)$ consists of all combinatorial surjective mappings and is generated by $\overline{D}_n$ together with the elementary regressions $\gamma_i{\times}\mathbf{Z}^{n-1}, ... \mathbf{Z}^{n-1}{\times}\gamma_i$.

The endofunctor $\mathcal{P}$ of *strongly reduced paths* is then defined as a colimit based on the category $\mathbf{c} \subset \mathbf{Tol}^{op}$, where $\mathbf{c}^{op} \supset \mathbf{d}^{op}$ consists of all *combinatorial surjective mappings* c: [i', j'] → [i, j]

(3) $\mathcal{P}X = \mathrm{Colim}_{\mathbf{c}}\, X^{[i, j]}$, $\qquad |\mathcal{P}X| = \mathrm{Colim}_{\mathbf{c}}\, !\mathbf{Smp}([i, j], X) = |PX|/\Gamma_1$.

Plainly, $|\mathcal{P}X|$ is a groupoid.

**3.6. Theorem (Preservation properties).** *The endofunctors* **P** *and* $\mathcal{P}$ *preserve equalisers, but do not preserve binary products nor pullbacks. They commute with exponentials* $(–)^{\hat{c}}$, *but not with* $(–)^{\mathbf{S}^0}$ *or* $(–)^{[h, k]}$ (*for* $k - h \geq 2$); *they are not special, nor pseudo-special* (1.6).

**Proof.** We consider **P**; the argument for $\mathcal{P}$ is analogous. The fact that **P**: !**Smp** → !**Smp** preserves equalisers need only be verified on |**P**| (1.5). For a pair of maps f, g: X → Y and a path $a \in |PX|$, the condition $fa \equiv ga$ means that fa.d = ga.d', for suitable delay functions d, d': **Z** → **Z**; we can thus replace a with a representative a' = ade = ad'e', so that fa' = ga'. Thus, the equaliser of |**P**f|, |**P**g|: |**P**X| → |**P**Y| is |**P**E|, where E is the equaliser of f, g.

As to products, let $Y = X_1{\times}X_2$; trivially, the congruence ≡ in |PY| implies the product of the congruences in $|PX_k|$: if $a_k.d = b_k.d'$ for all k, then $a_k \equiv b_k$; but the converse is false: take two different paths $a \equiv a'$: **Z** → $X = X_1 = X_2$ and consider the paths <a, a>, <a, a'>: **Z** → $X^2$; they have congruent projections on X, but are not congruent: <a, a>.d = <a, a'>.e would give ad = ae, ad = a'e, whence ae = a'e and a = a', as all delays are surjective mappings. Since **P** preserves the terminal, pullbacks are not preserved either. We also know that $(–)^{\mathbf{S}^0} = (–)^{[1]+[1]} = (–)^2$ does not commute with our colimit, while $(–)^{\hat{c}}$ does, because our site has products (1.5).

Finally, for $(–)^{[0, 2]}$, the comparison mapping

(1) $\mathrm{Colim}_{\mathbf{d}}([0, 2], X^{[i, j]}) \to ([0,2], \mathrm{Colim}_{\mathbf{d}}\, X^{[i, j]})$

(is surjective but) need not be injective, as shown by this counterexample. Let X be a simplicial complex having three distinct points x, y, z with x!y and y!z (e.g., X = Δ[2] or [0, 2]); then the following two maps in X, defined on [0, 1]×[0, 2] and [0, 2]×[0, 2], respectively

$$(1)\qquad \begin{matrix} x & y \\ y & y \\ y & z \end{matrix} \qquad\qquad \begin{matrix} x & \mathbf{x} & y \\ y & \mathbf{y} & \mathbf{y} \\ y & z & \mathbf{z} \end{matrix} \qquad\qquad \begin{matrix} 0 \\ 1 \\ 2 \end{matrix}$$

represent different elements of $([0, 2], \mathrm{Colim}_{\mathbf{d}} X^{[i, j]})$, since the vertical triple $(x, y, z)$ only appears in the second; but they are identified in $\mathrm{Colim}_{\mathbf{d}}([0, 2], X^{[i, j]})$, where we can convert the first map into the second by two delays, on $[0, 1]$ and $[1, 2]$ (the 'delayed parts' are in boldface characters).

**3.7. The fundamental ω-groupoid.** Consider now the set-valued functors $|P^nX|/{\equiv_n}$ produced by the congruence defined by the submonoid $D_n \subset \mathbf{Cs}(\mathbf{Z}^n, \mathbf{Z}^n)$ of n-*dimensional delays*. $D_2$ is generated by the elementary delays $\delta_i{\times}\mathbf{Z}$, $\mathbf{Z}{\times}\delta_{i'}$ $(i, i' \in \mathbf{Z})$ and still satisfies the cofiltering property 2.3.2 (the generators of different families commute); and so on.

Examining the last counterexample (3.6.1), we see that $|\mathbf{P}^2(X)|$ is different from $|P^2X|/{\equiv_2}$: the latter is a symmetric construct, with a symmetry induced by $s\colon P^2 \to P^2$, the former is not (it is a quotient of $|P^2X|$ by a stronger equivalence relation, *non-invariant under* $s$, where 'waving delays' are permitted in *one* direction). Therefore, we cannot transfer to $\mathbf{P}$ the previous structure of $P$ and we shall 'replace' its powers $\mathbf{P}^nX$ with the *set-valued* functors

(1) $\mathbf{P}^{(n)}X = \mathrm{colim}_{\mathbf{d}^n}\, !\mathbf{Smp}([i_1, j_1] \times ... \times [i_n, j_n], X) = |P^nX|/{\equiv_n}$ ,

produced by the congruences $\equiv_n$. Of course, $\mathbf{P}^{(0)}(X) = |X|$ and $\mathbf{P}^{(1)}(X) = |\mathbf{P}X|$.

For each object $X$, these functors produce a *cubical set* $\mathbf{P}^*X$ *with connections* $(g_i)$ [4, 1] *and symmetries* (the interchanges $s_i$ and the reversions $r_i$)

(2) $|X| = \mathbf{P}^{(0)}(X) \;\leftleftarrows\; \mathbf{P}^{(1)}(X) \;\overset{\textstyle\leftleftarrows}{\leftleftarrows}\; \mathbf{P}^{(2)}(X) \;...\; \leftleftarrows\; \mathbf{P}^{(n)}(X) \;...$

(only faces are shown); the structural maps (for $1 \le i \le n$; $\kappa = \pm$)

(3) $\partial^\kappa_i\colon \mathbf{P}^{(n)} \to \mathbf{P}^{(n-1)}$, $\qquad e_i\colon \mathbf{P}^{(n-1)} \to \mathbf{P}^{(n)}$, $\qquad g^\kappa_i\colon \mathbf{P}^{(n)} \to \mathbf{P}^{(n+1)}$,
$\quad s_i\colon \mathbf{P}^{(n+1)} \to \mathbf{P}^{(n+1)}$, $\qquad r_i\colon \mathbf{P}^{(n)} \to \mathbf{P}^{(n)}$,

are induced by five natural transformations of the powers of $P$

(4) $\partial^\kappa_i = P^{n-i}\partial^\kappa P^{i-1}\colon P^n \to P^{n-1}$, $\qquad e_i = P^{n-i}eP^{i-1}\colon P^{n-1} \to P^n$,
$\quad g^\kappa_i = P^{n-i}g^\kappa P^{i-1}\colon P^n \to P^{n+1}$,
$\quad s_i = P^{n-i}sP^{i-1}\colon P^{n+1} \to P^{n+1}$, $\qquad r_i = P^{n-i}rP^{i-1}\colon P^n \to P^n$.

Moreover, there is an i-*composition* $k_i\colon \mathbf{P}^{(n)} \times_{\mathbf{P}^{(n-1)}} \mathbf{P}^{(n)} \to \mathbf{P}^{(n)}$ (pullback with respect to the i-faces $\partial^\kappa_i\colon \mathbf{P}^{(n)} \to \mathbf{P}^{(n-1)}$), induced by the concatenation in direction i of the functor $\mathbb{P}$ of Moore paths, as for $n = 1$, in 3.4 (which is the main reason why we need here to consider $\mathbb{P}$)

(5) $k_i = \mathbb{P}^{n-i}k\mathbb{P}^{i-1}: \mathbb{P}^{n-i}\mathbb{Q}\mathbb{P}^{i-1} \to \mathbb{P}^n$.

Finally, $\mathbf{P}^*X$ is a (small) *cubical* ω*-category with connections* [1] *and symmetries*. It contains an *involutive ω-category* $\mathbf{P}_*X$ (with involutions $r_i$ for all composition laws)

(6) $|X| = \mathbf{P}_0(X) \;\leftrightarrows\; \mathbf{P}_1(X) \;...\; \leftrightarrows\; \mathbf{P}_n(X) \;...$

obtained in the obvious way, extending 2.4.2. First $\mathbf{P}_i(X) = \mathbf{P}^{(i)}(X)$ for $i = 0, 1$; then $\mathbf{P}_2(X) \subset \mathbf{P}^{(2)}(X)$ is the subset of those double reduced paths whose first faces are degenerate; $\mathbf{P}_3(X) \subset \mathbf{P}^{(3)}(X)$ contains those triple cells whose last faces end in $\mathbf{P}_2(X)$, while the first and second are degenerate in their last direction; and so on (always keeping in the structure the last faces):

(7) $a \in \mathbf{P}_n(X)$ iff: - $a \in \mathbf{P}^{(n)}X$, $\quad \partial^\varepsilon(a) = \partial^\varepsilon_n(a) \in \mathbf{P}_{n-1}X$,
- $\partial^\varepsilon_i(a) \in \mathrm{Im}(e^{n-1}_{n-1}: \mathbf{P}^{(n-2)}X \to \mathbf{P}^{(n-1)}X)$, for $0 \leq i < n$.

Each component has involutions $r_i: \mathbf{P}^{(n)} \to \mathbf{P}^{(n)}$, for $1 \leq i \leq n$ (which can be extended, letting $r_i = 1$ for $i > n$). Note that $\mathbf{P}^*X$ and $\mathbf{P}_*X$ 'contain the same information': it has been quite recently proved by Al-Agl - Brown - Steiner [1] that (edge-symmetric) *cubical* ω*-categories with connections* are equivalent to ω*-categories* (by the procedure described in (7)).

Replacing $\mathbf{P}$ with $\mathcal{P}$ (and $\bar{D}_n$ with $\Gamma_n$), in these constructs, we obtain a (strict) *cubical* ω*-groupoid with connections and interchanges* $\mathcal{P}^*X$

(8) $\mathcal{P}^{(n)}X = \mathrm{colim}_{\mathbf{c}^n}\, !\mathbf{Smp}([i_1, j_1] \times ... \times [i_n, j_n], X) = |P^nX|/\Gamma_n$,

and, within the latter, the *fundamental* ω*-groupoid* $\Pi_\omega X = \mathcal{P}_*X$, defined as in (7). Again, $\mathcal{P}^*X$ and $\mathcal{P}_*X$ 'contain the same information', since cubical ω-groupoids with connections are equivalent to (globular) ω-groupoids, a fact already known from Brown - Higgins' results [4, 6].

**3.8. Fundamental n-groupoids.** Finally, the *fundamental* n*-groupoid* of a symmetric simplicial set X is produced by the reflector $\rho_n$: ω-**Gpd** → n-**Gpd** of the (skeletal) embedding of n-**Gpd**

(1) $\Pi_n(X) = \rho_n\Pi_\omega(X)$.

Explicitly, $\Pi_n(X)$ only differs in degree n from the n-globular set produced by truncation of $\mathcal{P}_*X$

(2) $|X| = \mathcal{P}_0(X) \;\leftleftarrows\!\!\rightarrow\; \mathcal{P}_1(X) \;\dots\; \leftleftarrows\!\!\rightarrow\; \mathcal{P}_{n-1}(X) \;\leftleftarrows\!\!\rightarrow\; \Pi_n(X)$

and the *set* $\Pi_n(X)$ (containing all lower cells as degenerate elements) is the coequaliser of the faces ending in $\mathcal{P}_n(X)$, i.e. the quotient $\mathcal{P}_n(X)/\simeq_{n+1}$ modulo the equivalence relation $a \simeq_{n+1} b$: there exists $A \in \mathcal{P}_{n+1}(X)$ with $\partial^- A = a$, $\partial^+ A = b$. The fundamental groupoid $\Pi_1$ coincides with the structure $|PX|/\simeq_2$ already introduced in the previous Section (2.5), because $\Gamma_1$-congruence implies the 2-homotopy relation used there (3.5).

For a *pointed* symmetric simplicial set X, the homotopy group $(n \geq 1)$

(3) $\pi_n(X) = \Pi_n(X)(*_X, *_X) = \Omega^{(n)}(X)/\simeq_{n+1}$,

is the group of endocells of $\Pi_n(X)$ at the degenerate n-tuple strongly reduced path at the base point. By the usual argument, it is abelian for $n \geq 2$ (when all concatenation laws coincide). It can be viewed as the quotient of the involutive n-category $\Omega^{(n)}(X)$ of n-loops (those n-tuple strongly reduced paths $a \in \mathcal{P}^{(n)}(X)$ whose faces are constant at the base point), modulo $\simeq_{n+1}$.

Similarly (but this construct will not be used here), one can consider a *weak fundamental* n-*groupoid* $W\Pi_n(X) = \rho_n \mathbf{P}_*(X)$ defined as the involutive n-category produced by the reflector $\rho_n$: ω-**Cat** → n-!**Cat** of the embedding of involutive n-categories. Actually, $W\Pi_n(X)$ is a sort of *semi-strict* n-groupoid: using connections and reversions, one can see that

- $W\Pi_n(X)$ is a *strict* n-involutive n-category;

- the last involution $r_n$ provides n-cells with *strict* inverses with respect to the last composition $k_n$;

- $r_{n-1}$ provides (n–1)-cells with quasi-inverses, up to assigned (invertible!) n-cells; and so on.

**3.9. Pasting Theorem.** *The functors*

(1) $\mathcal{P}^*$: !**Smp** → ω-!C**Gpd**, $\Pi_\omega = \mathcal{P}_*$: !**Smp** → ω-**Gpd**,
$\Pi_n = \rho_n.\mathcal{P}_*$: !**Smp** → n-**Gpd**,

*preserve all colimits* (ω-!C**Gpd** *denotes the category of cubical ω-groupoids with connections and symmetries*). *They have right adjoints. In particular, the* symmetric ω-nerve, *right adjoint to* $\Pi_\omega$

(2) $M_\omega$: ω-**Gpd** → !**Smp**, $M_\omega(A)$ = ω-**Gpd**$(\Pi_\omega !\Delta, A)$: $!\mathbb{A}^{op}$ → **Set**,

*is 'represented' by the restriction* $\Pi_\omega.!\Delta\colon !\mathbb{\Delta} \to \omega$-**Gpd** *of* $\Pi_\omega$ *to the Yoneda embedding. By composition,* $\Pi_n$ *is left adjoint to the* symmetric n-nerve, *restriction of* $M_\omega$

(3) $\Pi_n\colon$ !**Smp** $\rightleftarrows$ n-**Gpd**: $M_n$, $\qquad\qquad\qquad \Pi_n \dashv M_n$.

**Proof.** As to the first assertion, since sums are plainly preserved and $\rho_n$ is a left adjoint, we only need to prove that $\mathcal{P}^*$ and $\Pi_\omega = \mathcal{P}_*$ preserves coequalisers. It is sufficient to do so for the former, because of the already recalled Brown-Higgins equivalence (3.7).

Let f, g: X $\to$ Y have coequaliser p: Y $\to$ Z = Y/R in !**Smp** (R is the least congruence of symmetric simplicial sets such that pf = pg), preserved by any truncation $tr_n$: **Smp** $\to$ $\mathbf{Smp}_n$ (a left adjoint). Now, let h: $\mathcal{P}^*Y \to C$ coequalise $\mathcal{P}^*f$, $\mathcal{P}^*g$ in $\omega$-!C**Gpd**, and let us prove that it factorises through the $\omega$-functor $\mathcal{P}^*p\colon \mathcal{P}^*Y \to \mathcal{P}^*Z$, by a unique $\omega$-functor k: $\mathcal{P}^*Z \to C$.

At the level of vertices, $h_0\colon Y_0 \to C_0$ coequalises $f_0$, $g_0$, and factorises uniquely as $k_0p_0\colon Y_0 \to Z_0 \to C_0$ by the previous remark on truncation. At the level of 1-cells, a reduced directed path $z^\bullet \in \mathcal{P}_1Z$ (z: **Z** $\to$ Z) is represented by a finite sequence $z = (z_i)$ of consecutive 1-simplices of Z ($\partial_0 z_{i-1} = \partial_1 z_i$); since Z is a quotient of Y, the sequence has a lifting $y = (y_i)$ in $Y_1$ ($p(y_i) = z_i$); now, the lifted sequence *need not be a path*, but certainly its image $(h_1(y_i))$ via $h_1\colon \mathcal{P}^{(1)}Y \to C_1$ is a sequence of consecutive arrows, whose composite $k_1(z^\bullet)$ does not depend on the choice of z and y: we have thus a unique mapping $k_1\colon \mathcal{P}^{(1)}Z \to C_1$ such that $k_1p_1 = h_1$. For the 2-dimensional level, a reduced directed double path $z^\bullet \in \mathcal{P}^{(2)}Z$ is represented by a double path z: $\mathbf{Z}^2 \to Z$; this is a diagram of 2-simplices of Z, parametrised over the 2-simplices of $\mathbf{Z}^2$ contained in a finite rectangle

(4) 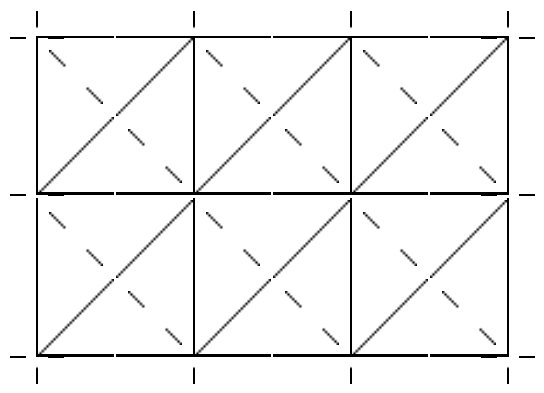

coherently with faces and symmetries. Again, a lifting of z in $Y_2$ need not be coherent (a double path) but, transformed by $h_2\colon \mathbf{P}^{(2)}Y \to C_2$, becomes a double path, and $k_2(z^\bullet)$ is defined as its pasting in the involutive double category $\rho_2C$. And so on.

Now, the existence of the right adjoint $M_\omega$ follows from general facts about an arbitrary presheaf category $\mathrm{Psh}\mathcal{C} = \mathbf{Set}^{\mathcal{C}^{op}}$, well known since the introduction of adjoint functors by Kan [17, 18] (see also [22], I.5). The classical starting point is the topological model of ordinals $\Delta: \mathbb{\Delta} \to \mathbf{Top}$, which produces the singular simplicial functor $S_* = \mathbf{Top}(\Delta, -): \mathbf{Top} \to \mathrm{Psh}\mathbb{\Delta}$ and, as a left adjoint to the latter, the geometric realisation of simplicial sets. For the general Kan's result, to be applied here, replace $\Delta$ with any functor $\Sigma: \mathcal{C} \to \mathbf{T}$ from a small category to a category *with (small) colimits*. Then, the functor 'represented' by $\Sigma$

(5) $H^\Sigma: \mathbf{T} \to \mathbf{Set}^{\mathcal{C}^{op}}$, $\qquad H^\Sigma(T) = \mathbf{T}(\Sigma(-), T): \mathcal{C}^{op} \to \mathbf{Set}$,

has a left adjoint, which can be calculated as a coend over $\mathcal{C}$

(6) $\otimes\Sigma: \mathbf{Set}^{\mathcal{C}^{op}} \to \mathbf{T}$, $\qquad \otimes\Sigma(X) = \int^c X(c){\cdot}\Sigma(c)$.

(Or as a colimit over the category of elements of $X$, as in [22], I.5. $\otimes\Sigma$ is the original notation in [18], Section 3. If $\Sigma = y: \mathcal{C} \to \mathrm{Psh}\mathcal{C}$, then $H^\Sigma$ is the identity, whence also $\otimes\Sigma$ is so, and $X = \int^c X(c){\cdot}y(c)$ is the usual expression of an arbitrary presheaf as a colimit of representable ones.)

In other words, a functor $F: \mathbf{Set}^{\mathcal{C}^{op}} \to \mathbf{T}$ with values in a cocomplete category has a right adjoint iff it *preserves colimits,* and then $F \cong \otimes\Sigma \dashv H^\Sigma$, where $\Sigma = Fy: \mathcal{C} \to \mathbf{T}$

(7) $F(X) \cong F(\int^c X(c){\cdot}y(c)) \cong \int^c X(c){\cdot}Fy(c) = \otimes\Sigma(X)$,

so that $\otimes\Sigma$ is characterised, up to isomorphism, as the colimit-preserving functor which extends $\Sigma$ along the Yoneda embedding, $\Sigma \cong \otimes\Sigma.y$.

## 4. Directed combinatorial homotopy, in degree 1

The *directed line* $\uparrow\mathbf{Z}$ gives a notion of (non-reversible) *directed combinatorial homotopy* in **Smp**.

**4.1. Directed combinatorial spaces.** The simple presheaves of **Smp** and $\mathbf{Smp}_1$ (1.3) have descriptions similar to their symmetric counterpart (1.4).

A *directed simplicial complex*, or *directed combinatorial space* (d-space for short), will be a set $X$ equipped with a set $\uparrow X$ of words in $X$, $x = (x_1,... x_p)$, called *linked words*, which: contains the empty sequence, contains all unary words (identified with points), and is closed under omitting or repeating entries. A *map* in their category $\uparrow\mathbf{Cs}$ is a mapping which preserves linked words.

$\uparrow\mathbf{Cs}$ has all limits and colimits, and is cartesian closed. $\mathrm{Hom}(X, Y) = Y^X$ is the set of maps, with the following directed combinatorial structure:

(1) a word $f$ in $\uparrow$**Cs**$(X, Y)$ is linked if, for every word $g = (g_1,... g_p)$ obtained by (possibly) repeating terms of $f$ and every linked word $x = (x_1,... x_p) \in \uparrow X$ (of the same length), the word $g(x) = (g_1(x_1),... g_p(x_p))$ is linked in $Y$.

Again, $\uparrow$**Cs** has a canonical cc-embedding in **Smp**, produced by the obvious realisation $o: \mathbb{A} \to \uparrow$**Cs** of ordinals as d-spaces (the linked words are the increasing ones), which identifies a directed combinatorial space with a *simple* simplicial set

(2) $c: \uparrow\mathbf{Cs} \to \mathbf{Smp}, \qquad (cA)_n = \mathbf{Smp}(o[n], A) = \{\text{linked n-words in } A\},$

with reflector (which preserves finite products)

(3) $u: \mathbf{Smp} \to \uparrow\mathbf{Cs}, \qquad uX = X_0 \ (\mathbf{x} \text{ linked iff } \exists\, x \in X_n\text{: } \mathbf{x} = (i^*x)_{0 \le i \le n}).$

Finally, a *step set* has a *step* relation or precedence relation $\prec$ (only assumed to be reflexive), and amounts to a simple graph; they form the category **Stp**, again complete, cocomplete and cartesian closed, with $f \prec g$ whenever $x \prec x'$ implies $f(x) \prec g(x')$. The diagram 1.4.3 is replaced with the following commutative square of cc-embeddings, having a similar description

(4)
$$\begin{array}{ccc} \mathbf{Stp} & \overset{s}{\underset{c}{\leftrightarrows}} & \uparrow\mathbf{Cs} \\ c \downarrow \uparrow u & & c \downarrow \uparrow u \\ \mathbf{Smp}_1 & \overset{tr}{\underset{cosk}{\leftrightarrows}} & \mathbf{Smp} \end{array}$$

**4.2. Directed line and paths.** Directed combinatorial homotopy, or d-homotopy, will be defined in each of these categories (**Smp**, $\mathbf{Smp}_1$, $\uparrow$**Cs**, **Stp**) by means of the *directed integral line* $\uparrow\mathbf{Z}$, the set of integers with the structure of *consecutivity*, and its powers.

To begin with, $\uparrow\mathbf{Z}$ is a step set, with $i \prec j$ if $0 \le j - i \le 1$; as a graph, $\uparrow\mathbf{Z}$ has precisely one map $i \to j$ whenever $i \prec j$; as a directed simplicial complex, a word $(i_1,... i_n)$ is linked if $i_p \prec i_q$ for all $p \le q$, iff it is of type $(i,... i, i+1,... i+1)$ or constant; finally, as a simplicial set, its simplices are sequences of the same type. Again, a finite integral interval $\uparrow[i, j] \subset \uparrow\mathbf{Z}$ is finitely presentable. The elementary directed interval $\uparrow[0, 1]$ is the step-subset $\{0, 1\}$, i.e. the ordinal $[1]$ and the representable presheaf $\Delta[1]$. The lattice operations $\vee, \wedge: (\uparrow\mathbf{Z})^2 \to \uparrow\mathbf{Z}$ preserve the step relation; the same holds for translations.

Of course, $\uparrow\mathbf{Z}$ is not isomorphic to the integral line with the natural order, which is why we prefer the term 'directed' to 'ordered'. 'Zigzag' step sets show that also the term 'oriented d-space' could be misleading.

We have thus the representable functor of *directed lines*

(1) $\uparrow L = (-)^{\uparrow\mathbf{Z}}: \mathbf{Smp} \to \mathbf{Smp}, \qquad (\uparrow LX)_n = (X^{\uparrow\mathbf{Z}})_n = \mathbf{Smp}(\uparrow\mathbf{Z}, X^{\Delta[n]}).$

In its underlying set, much as is the symmetric case (2.2.2)

(2) $|\uparrow LX| = \mathbf{Smp}(\uparrow\mathbf{Z}, X) = \mathbf{Smp}(sk_1tr_1\uparrow\mathbf{Z}, X) = \mathbf{Smp}_1(tr_1\uparrow\mathbf{Z}, tr_1X),$

a (directed) line $a: \uparrow\mathbf{Z} \to X$ amounts to its 1-truncation $(a_i: \uparrow\mathbf{Z}_i \to X_i)_{i\leq 1}$, a morphism of graphs, i.e. a sequence of consecutive 1-simplices $a_1(i, i+1) \in X_1$; again, this line is said to be a *path* if there exists a support $\rho = [\rho^-, \rho^+] \subset \mathbf{Z}$ such that the 1-simplex $a_1(i, i+1)$ is degenerate for $i < \rho^-$ and for $i \geq \rho^+$. Paths form a subfunctor $\uparrow P_\bullet \subset \uparrow L_\bullet: \mathbf{Smp} \to \mathbf{Set}$, which lifts canonically to a *directed-path endofunctor* $\uparrow P: \mathbf{Smp} \to \mathbf{Smp}$, a special subfunctor of $\uparrow L$ (1.5-6)

(3) $\uparrow PX = \mathrm{Colim}_{\mathbf{t}}\, X^{\uparrow[i, j]}, \qquad |\uparrow PX| = \mathrm{Colim}_{\mathbf{t}}\, \mathbf{Smp}(\uparrow[i, j], X),$

since all truncations in $\mathbf{t}^{op}$ (3.1) are *increasing* functions: $\mathbf{t}^{op} \subset \mathbf{Stp} \subset \uparrow\mathbf{Cs} \subset \mathbf{Smp}$.

Line-functors $\uparrow LX = X^{\uparrow\mathbf{Z}}$ and paths are defined in the same way for the cc-subcategories **Stp**, $\uparrow$**Cs**, $\mathbf{Smp}_1$ (4.1.4). For d-spaces, a line $a: \uparrow\mathbf{Z} \to X$ is a mapping preserving the linked words, and a word $(a_1,... a_n)$ of lines is linked iff each word $(a_1(i),... a_p(i), a_p(i+1),... a_n(i+1))$ is linked in X. $|\uparrow PX|$ is the subobject of the lines $a: \uparrow\mathbf{Z} \to X$ eventually constant at the left and the right.

Directed paths in a *symmetric* simplicial set viewed in **Smp** are the same as symmetric paths: $\mathbf{Smp}(\uparrow[i, j], UX) = !\mathbf{Smp}(\mathrm{Sym}(\uparrow[i, j]), X) = !\mathbf{Smp}([i, j], X)$. (Here, $\mathrm{Sym} \dashv U$ and $\mathrm{Sym}(X) = \mathrm{Lan}(X)$ is the left Kan extension of the presheaf $X: \mathbb{A}^{op} \to \mathbf{Set}$ along $\mathbb{A}^{op} \subset !\mathbb{A}^{op}$).

**4.3. Directed homotopies.** A (directed) *homotopy* $a: f \to g: X \to Y$ *of simplicial sets* will be a map $a: X \to \uparrow PY$ with $\partial^- a = f$, $\partial^+ a = g$.

Our homotopy is said to be *immediate* if it factorises through $Y^{\uparrow[0, 1]}$, so that $a: X \to Y^{\uparrow[0, 1]}$ can also be viewed as $a: \uparrow[0, 1] \to Y^X$, an immediate path in $Y^X$. *Immediate homotopies amount to the classical notion of simplicial homotopy* [24], as it follows from the fact that $\uparrow[0, 1]\times X = \Delta[1]\times X$ is generated by the pairs $(e^n\varepsilon, e_ix) \in X_{n+1}$ $(x\in X_n,\ 0 \leq i \leq n)$. More generally, the homotopy is *bounded* if it factorises through some $Y^{\uparrow[i, j]}$, i.e. if it is a path in $Y^X$

(1) $|\uparrow P(Y^X)| = \mathrm{Colim}_{\mathbf{t}}\, \mathbf{Smp}(\uparrow[i, j], Y^X) = \mathrm{Colim}_{\mathbf{t}}\, \mathbf{Smp}(X, Y^{\uparrow[i, j]});$

note that, *if* X *is finitely generated* (= finitely presentable), this is the same as $|(\uparrow PY)^X|$, and all homotopies from X are bounded.

In the set **Smp**(X, Y) we have various 'directed homotopy' relations, defined by the existence of an *immediate* (or *bounded*, *bounded on connected components, arbitrary*) homotopy $f \to g$

(2) $f \prec_i g \;\Rightarrow\; f \prec_b g \;\Rightarrow\; f \prec_c g \;\Rightarrow\; f \prec_1 g.$

Generally, these relations are just reflexive. Nevertheless, bounded homotopies – as paths in $Y^X$ – can be concatenated and obtained as a finite concatenation of immediate ones: therefore, the relation $\prec_b$ *is* transitive and coincides with the preorder relation spanned by $\prec_i$. The equivalence relations generated by the previous reflexive relations will be written

(3) $f \simeq_b g \;\Rightarrow\; f \simeq_c g \;\Rightarrow\; f \simeq_1 g,$

but one should be aware that these relations are generally *rather weak*, as will be discussed at the end of the next section. For a *finitely generated* X, we already know that $\prec_b$, $\prec_c$ and $\prec_1$ coincide in any set **Smp**(X, Y); their equivalence relations coincide as well. On the other hand, *within* Kan complexes, it is well known that $\prec_i$ is an equivalence relation: it coincides thus with $\prec_b$ and $\simeq_b$.

**4.4. The fundamental category.** A directed combinatorial space X has a *fundamental category* $\uparrow\Pi(X)$, whose objects are the points of X and whose arrows $[a]: x \to x'$ are the relative-homotopy classes of paths, up to the equivalence relation $\simeq_2$. The latter is spanned by the relation $a' \prec_2 a''$: there exists a double directed path $A \in |\uparrow P^2X| \subset \mathbf{Smp}((\uparrow\mathbf{Z})^2, X)$ whose last faces are a', a'' while the first faces $\partial_1^\kappa A = (\uparrow P\partial^\kappa)A$ are degenerate.

This fundamental-category functor can be obtained as a left adjoint

(1) $\uparrow\Pi_1 : \mathbf{Smp} \rightleftarrows \mathbf{Cat} : N_1, \qquad\qquad \uparrow\Pi_1 \dashv N_1,$

where $N_1C$ is the *nerve* of a small category, i.e. the representable functor produced by the embedding $i: \mathbb{A} \to \mathbf{Cat}$ which sends the ordinal $[n] > 0$ to the corresponding category

(2) $N_1C = \mathbf{Cat}(i(–), C): \mathbb{A}^{op} \to \mathbf{Cat},$

$(N_1C)_n = \mathbf{Cat}(i[n], C)$

$= \{(a_{ij})_{0\leq i\leq j\leq n} \mid a_{ij} \in C_1;\; a_{ii} \text{ is an identity};\; a_{ij} + a_{jk} = a_{ik}\}.$

Again (as $\Pi_1$ in 2.6), our functor $\uparrow\Pi_1$ preserves all colimits; this general result restricts to a van Kampen theorem on $\uparrow$**Cs**, for a union of subobjects (a pushout over their intersection). Similarly, the functor $\uparrow\Pi_0$: **Smp** → **Set** of path components is left adjoint to the discrete embedding D. The higher homotopy functors $\uparrow\Pi_n$ will be dealt with in Section 5.

As in the symmetric case (2.1), there is in **Smp** and $\mathbf{Smp}_1$ a k-*point directed circle* $\uparrow C_k$ ($k \geq 1$), the coequaliser of the faces of $\uparrow[0, k]$; it belongs to **Stp** (and $\uparrow$**Cs**) for $k \geq 3$; moreover, $\mathrm{Sym}(\uparrow C_k) = C_k$ (4.2). The *directed 0-sphere* is symmetric, the discrete $\mathbf{S}^0 = \{-1, 1\} \subset \mathbf{Z}$ (pointed at 1). Again, there is no standard n-sphere in $\uparrow$**Cs**, for $n > 0$, but we can surrogate it by a system of directed combinatorial spaces, similar to the symmetric $C_kS^n$ (2.7): the *directed* k-*collapsed* n-*sphere* $\uparrow C_kS^n = \uparrow\mathbf{Z}^n/\sim_k$ is obtained by collapsing to a point all the points of $\uparrow\mathbf{Z}^n$ out of the cube $[1, k]^n$ (for $k \geq 2$), and can be used to present the homotopy monoid $\uparrow\pi_n(X)$ (5.1.6) as a colimit of the system $[\uparrow C_kS^n, X]$.

A directed homotopy $a: f \to g: X \to Y$ plainly gives $\uparrow\Pi_0(f) = \uparrow\Pi_0(g)$; but, at the level of fundamental categories, we just get a natural transformation $\uparrow\Pi_1(a): \uparrow\Pi_1(f) \to \uparrow\Pi_1(g)$, *which need not be invertible*. For instance, the endomaps 0, $\mathrm{id}: \uparrow\mathbf{N} \to \uparrow\mathbf{N}$ of the directed integral half-line are linked by a *telescopic* directed homotopy $a: 0 \to \mathrm{id}$, $a(i, t) = 0 \vee (i \wedge t)$ (cf. [11], 3.4), but do not induce isomorphic functors, via $\uparrow\Pi_1$. A study of such facts cannot be given here; we only hint at a stronger notion, sufficient to make $\uparrow\Pi_1(a)$ a natural equivalence: there exist two *bounded* homotopies a, b such that

(4) $\quad a: f \rightleftarrows g :b, \qquad\qquad a + b \simeq_2 0_f, \qquad b + a \simeq_2 0_g,$

where $\simeq_2$ denotes, as above but more generally, the equivalence relation spanned by $a' \prec_2 a''$ (there exists a map $A: X \to \uparrow P^2Y$ whose last faces are a', a'' while the first faces are degenerate). Higher homotopy functors require higher conditions.

**4.5. Step metric spaces.** Applications of the (reversible) homotopy theory of **Cs** to the analysis of images have been studied in [11]. A metric space X has a combinatorial structure $t_\varepsilon X$ *at each resolution* $\varepsilon \in [0, \infty]$, defined by the tolerance relation $d(x, y) \leq \varepsilon$; the resulting n-homotopy group $\pi_n^\varepsilon(X) = \pi_n(t_\varepsilon X)$ detects singularities which can be captured by an n-dimensional grid, with edges bound by $\varepsilon$; this works equally well for continuous euclidean regions or discrete ones, as produced by a scanner.

A 'directed' counterpart of these applications has been briefly considered in [12], 5.7-8. A *step metric space* $X = (X, d, \prec)$ is a metric space equipped with a step relation. It has a family of canonical step-structures $s_\varepsilon X$ *at resolution* $\varepsilon$, where $x \prec_\varepsilon x'$ iff $x \prec x'$ and $d(x, x') \leq \varepsilon$. Thus, a path in $s_\varepsilon X$ is based on a sequence of points $(x_0,... x_p)$ with $x_{i-1} \prec_\varepsilon x_i$; the homotopy theory of directed combinatorial spaces gives a family of fundamental categories *at resolution* $\varepsilon$, $\uparrow\Pi_1^{\varepsilon}(X) = \uparrow\Pi_1(s_{\varepsilon}X)$, which can be studied much as in [11] and has similar applications. Loosely speaking, $\uparrow\Pi_1^{\varepsilon}(X)$ is based on directed paths as above; their deformations, formed of 2-dimensional words, are similarly discrete.

## 5. Higher fundamental categories and conclusions

On the presheaf category **Smp**, all higher homotopy functors preserve colimits and are left adjoints to higher nerves (5.1-2). We conclude with some conjectures on such nerves (5.3) and some general comments on defining homotopy in presheaf categories (5.4).

**5.1. Higher fundamental categories.** Let us move now to higher homotopies in **Smp**; these can be treated as in the symmetric case (Section 3), except of course that reversion, regressions (3.9) and all their consequences are missing. Let X be a simplicial set.

Also here, we need to consider the pseudo-special endofunctor $\uparrow\mathbb{P}$ of directed Moore paths and the functor $\uparrow\mathbf{P}$ of *directed reduced paths* (just a canonical lifting), related by natural transformations

(1) $\uparrow\mathbb{P}X \to \uparrow PX \to \uparrow\mathbf{P}X,$

$\uparrow\mathbb{P}X = \Sigma_{\mathbf{s}} X^{\uparrow[i, j]}, \qquad |\uparrow\mathbb{P}X| = \Sigma_{\mathbf{s}} (\uparrow[i, j], X) \subset \mathbf{s}\times|\uparrow PX|,$

$\uparrow\mathbf{P}X = \mathrm{Colim}_{\mathbf{d}} X^{\uparrow[i, j]}, \qquad |\uparrow\mathbf{P}X| = \mathrm{Colim}_{\mathbf{d}}\, \mathbf{Smp}(\uparrow[i, j], X) = |\uparrow PX|/\!\equiv,$

since $\mathbf{s} \subset \mathbf{t} \subset \mathbf{d} \subset \mathbf{Stp}^{\mathrm{op}} \subset \mathbf{Smp}^{\mathrm{op}}$: all generalised delays are *step-maps* $\uparrow\mathbf{Z} \to \uparrow\mathbf{Z}$ (whereas $\mathbf{c} \not\subset \mathbf{Stp}^{\mathrm{op}}$). Similarly, $D_n \subset \mathbf{Stp}(\uparrow\mathbf{Z}^n, \uparrow\mathbf{Z}^n)$, and we have set-valued functors corresponding to 3.7.1

(2) $\uparrow\mathbf{P}^{(n)}X = \mathrm{colim}_{\mathbf{d}^n}\, \mathbf{Smp}([i_1, j_1] \times ... \times [i_n, j_n], X) = |\uparrow P^n X|/\!\equiv_n .$

The family $\uparrow\mathbf{P}^*(X) = (\uparrow\mathbf{P}^{(n)}(X))_{n\geq 0}$ of all n-dimensional directed reduced path sets is a cubical set with connections and interchange (but no reversion). The concatenation of consecutive reduced paths makes $\uparrow\mathbf{P}^*(X)$ into a cubical ω-category with interchange (only faces and degeneracies are drawn)

(3) $|X| = \uparrow\mathbf{P}^{(0)}(X) \quad \leftleftarrows\!\!\rightarrow \quad \uparrow\mathbf{P}^{(1)}(X) \; ... \; \overset{\leftarrow}{\underset{\rightarrow}{\leftrightarrows}} \quad \uparrow\mathbf{P}^{(n)}(X) \quad \ldots$

This contains the *fundamental* ω*-category* $\uparrow\Pi_\omega(X) = \uparrow\mathbf{P}_*(X)$, defined as in 3.7.7

(4) $|X| = \uparrow\mathbf{P}_0(X) \quad \leftleftarrows\!\!\rightarrow \quad \uparrow\mathbf{P}_1(X) \quad ... \quad \leftleftarrows\!\!\rightarrow \quad \uparrow\mathbf{P}_n(X) \quad ...$

The reflector $\rho_n$: ω-**Cat** → n-**Cat** produces the *fundamental* n-*category* $\uparrow\Pi_n(X) = \rho_n.\uparrow\Pi_\omega(X)$

(5) $|X| = \uparrow\mathbf{P}_0(X) \leftrightarrows \uparrow\mathbf{P}_1(X) \ \ ... \ \leftrightarrows \uparrow\mathbf{P}_{n-1}(X) \leftrightarrows \uparrow\Pi_n(X);$

again, the *set* $\uparrow\Pi_n(X)$ is the coequaliser of the faces $\uparrow\mathbf{P}_{n+1}(X) \rightrightarrows \uparrow\mathbf{P}_n(X)$, i.e. the quotient $\uparrow\mathbf{P}_n(X)/\simeq_{n+1}$ modulo the equivalence relation spanned by the relation $a \rightarrow_{n+1} b$: there exists $A \in \uparrow\mathbf{P}_{n+1}(X)$ with $\partial^-A = a$, $\partial^+A = b$.

For a *pointed* simplicial set X, the n-dimensional *homotopy monoid* is the monoid of endomaps

(6) $\uparrow\pi_n(X) = \uparrow\Pi_n(X)(*_X, *_X) = \uparrow\Omega^{(n)}(X)/\simeq_{n+1}$ (n ≥ 1),

at the degenerate n-tuple reduced path at the base point; it is abelian for n ≥ 2.

**5.2. Pasting Theorem.** *The functors* $\uparrow\Pi_\omega = \uparrow\mathbf{P}_*$: **Smp** → ω-**Cat** *and* $\uparrow\Pi_n = \rho_n.\uparrow\Pi_\omega$: **Smp** → n-**Cat** *preserve all colimits.* $\uparrow\Pi_\omega$ *has a right adjoint* $N_\omega$, *the* ω-nerve functor

(1) $N_\omega$: ω-**Cat** → **Smp**, $N_\omega(A)$ = ω-**Cat**($\uparrow\Pi_\omega.\Delta$, A): $\mathbb{\Delta}^{op}$ → **Set**,

*'represented' by the restriction* $\uparrow\Pi_\omega.\Delta$: $\mathbb{\Delta}$ → ω-**Cat** *of* $\uparrow\Pi_\omega$ *to the Yoneda embedding. The fundamental n-category functor is, by composition, left adjoint to the* n-nerve functor, *restriction of* $N_\omega$

(2) $\uparrow\Pi_n = \rho_n.\uparrow\mathbf{P}_*$: **Smp** ⇄ n-**Cat**: $N_n$, $\uparrow\Pi_n \dashv N_n$.

**Proof**. As for the symmetric case (3.9).

**5.3. Comments.** Again, the property of preserving colimits for $\uparrow\Pi_n$ is important in itself: it is a strong, simple version of the van Kampen theorem in **Smp**, and reduces to a standard version for $\uparrow$**Cs**, as in 2.6. Its consequence, the existence of the right adjoint, is perhaps even more relevant.

It would be interesting to determine *whether* $\uparrow\Pi_\omega.\Delta[n]$ *coincides with Street's* oriental $\mathcal{O}_n$, and more generally, *the functor* $\uparrow\Pi_\omega.\Delta$: $\mathbb{\Delta}$ → ω-**Cat** *coincides with the cosimplicial object* $\mathcal{O}$ *of orientals* ([27], Section 5). As a consequence, $N_\omega$ *would coincide with Street's ω-nerve* $N_\omega(A)$ = ω-**Cat**($\mathcal{O}$, A) *and* $\mathcal{O}_n$ *could also be viewed as the fundamental n-category* $\uparrow\Pi_n(\Delta[n])$ *of the n-simplex*. In fact, the latter has dimension n, so that the relation $\simeq_{n+1}$ is trivial in $\uparrow\mathbf{P}_n(X)$.

$\mathcal{O}_n$ is explicitly constructed in [27], and characterised as the free n-category generated by an n-graph formed by all subsets of [n] (p. 327). Another explicit construction can be found in [16].

**5.4. Homotopy for combinatorial sites.** Finally, some comments on how to construct 'combinatorial homotopy' in various presheaf categories, along the same lines.

Let us begin noting that three basic families of 'combinatorial (discrete) sites', the *simplicial*, *cubical* and *globular* family, present themselves naturally (in two versions, directed and symmetric) as inspired by simple geometric models: the tetrahedra $\Delta^n$, the cubes $\mathbf{I}^n = [0, 1]^n$, and the discs $\mathbf{D}^n$. In these cases, one can always realise the discrete site $\mathcal{C}$ as a subcategory of **Set** containing the terminal cardinal $[0] = 1$ and all mappings $c \to [0]$; the objects form a sequence of finite sets $[0], [1],...$ representing the geometric models of the corresponding topological dimension.

In the simplicial cases, already considered, our sequence consists of the finite positive cardinals $[n] = \{0, 1,... n\}$; but one might equivalently use for $[n]$ the integral trace of the tetrahedron $\Delta^n \cap \mathbf{Z}^{n+1} = \{e_0,... e_n\}$ (the unit points of the cartesian axes). In the cubical case, one can similarly use the *elementary cubical sets* $[n] = 2^n = \{0, 1\}^n = \mathbf{I}^n \cap \mathbf{Z}^n$, i.e. the integral traces of the standard cubes; the maps are conveniently defined, according to whether we are considering the directed or symmetric option. Finally, in the globular case, one can use the integral traces of the standard discs, $[n] = \mathbf{D}^n \cap \mathbf{Z}^n = \{\pm e_1,... \pm e_n\}$, and again choose the relevant maps as needed.

However, the basic facts – perhaps a hint for a definition of 'combinatorial site' – which have been used here to develop a reversible homotopy theory in **!Smp** (resp. a directed homotopy theory in **Smp**) seem to belong to a wide class of sites, satisfying two conditions (cf. 1.2, 1.6)

(A) in the corresponding presheaf topos, a finite product of representable presheaves is finitely presentable; moreover, the terminal presheaf is representable;

(B) the site contains a *standard reversible interval* (resp. a *standard directed interval*) [1], with the structure of a commutative, involutive cubical monoid *in* Psh$\mathcal{C}$

$$(1) \qquad y[0] \overset{\partial^\kappa}{\rightrightarrows} y[1] \overset{g^\kappa}{\leftleftarrows} y[1] \times y[1] \qquad\qquad r: [1] \to [1] \qquad (\kappa = \pm),$$

under axioms dual to 3.2.4 (resp.: omit the involution $r$); all representables are connected with respect to this interval (cf. 1.2 ($B_0$)).

(As already noted in 1.6, this implies that a finite product of representables is a *connected* finite colimit of representables.) On such bases, one constructs the system of 'finite intervals' as in 1.2, defines the path functor as a filtered colimit $PX = \mathrm{Colim}\, X^{[i, j]}$ (resp. $X^{\uparrow[i, j]}$) and the reduced path functor as in 3.4.1. A general, abstract construction might result quite complicated (in the extension of the algebraic structure); whilst, in concrete cases, one can take advantage of a simple description of the intervals, the line, truncations and delays.

Conditions (A), (B) are also satisfied, for instance, by the discrete site $\mathcal{O}$ of *finite connected ordered sets* (with underlying set a positive cardinal) and increasing mappings, which contains various models in each dimension $> 1$. Choosing the standard interval to be the ordinal 2, one would just get the previous theory applied to the underlying simplicial set; but other choices are possible.

Finally, the concrete relevance of these conditions in Homotopy is obvious for (B), but also clear for (A), which intervenes whenever a product $\Delta^p \times \Delta^q$ of standard simplices is covered with (p+q)-dimensional simplices intersecting at faces (a finite connected colimit of representable presheaves, viewed via the geometric realisation – a left adjoint).

Dipartimento di Matematica
Università di Genova
via Dodecaneso 35
16146 GENOVA, Italy
grandis@dima.unige.it